\documentclass[11pt,a4paper,twocolumn]{article}

\usepackage[sc]{mathpazo} % Use the Palatino font
\usepackage[T1]{fontenc} % Use 8-bit encoding that has 256 glyphs
\usepackage{microtype} % Slightly tweak font spacing for aesthetics

\newtheorem{theorem}{ Theorem}[section]
\newtheorem{definition}[theorem]{Definition}
\newtheorem{example}[theorem]{Example}
\usepackage[utf8]{inputenc}
\usepackage[english]{babel}
\usepackage{amsmath}
\usepackage{amsfonts}
\usepackage{amssymb}
\usepackage{makeidx}
\usepackage{graphicx}
\usepackage{lmodern}
\usepackage{kpfonts}
\usepackage[left=1.3cm,right=1.3cm,top=2cm,bottom=2cm]{geometry}
\usepackage{dsfont}
\usepackage{hyperref}
\usepackage{verbatim}
\usepackage{cite} 
\usepackage{caption}
\usepackage{subcaption}
\usepackage{fancyhdr}
\providecommand{\nset}[1]{
\mathbb{#1}
}
\providecommand{\set}[1]{
\left\{#1\right\}
}
\providecommand{\com}[1]{``#1"}

\providecommand{\ifr}[5]{
{}_{#1}^{#2}{#3}_{#4}^{#5}
}
\providecommand{\gam}[1]{
\Gamma\left(#1 \right)
}
\providecommand{\re}[1]{
Re\left(#1 \right)
}
\providecommand{\norm}[1]{
\left\lVert #1 \right\rVert
}

\providecommand{\ds}[1]{
\displaystyle{#1}
}

\usepackage{enumitem} % Customized lists
\setlist[itemize]{noitemsep} % Make itemize lists more compact

\usepackage{abstract} % Allows abstract customization
\usepackage{titlesec} % Allows customization of titles
\titleformat{\section}[block]{\large\scshape\centering}{\thesection.}{1em}{} % Change the look of the section titles
\titleformat{\subsection}[block]{\large\scshape\centering}{\thesubsection.}{1em}{} % Change the look of the section titles
\usepackage{titling} % Customizing the title section

\pretitle{\begin{center}\Huge\bfseries} % Article title formatting
\posttitle{\end{center}} % Article title closing formatting
\title{Proposal for Use the Fractional Derivative of Radial Functions in Interpolation Problems} % Article title
\author{
\textsc{A. Torres-Hernandez } \\ 
\normalsize Department of Physics - UNAM \\ 
\normalsize  \href{mailto:anthony.torres@ciencias.unam.mx}{anthony.torres@ciencias.unam.mx}
\and  
\textsc{F. Brambila-Paz }\\ 
\normalsize Department of Mathematics - UNAM \\
\normalsize \href{mailto:fernandobrambila@gmail.com}{fernandobrambila@gmail.com} 
\and  
\textsc{C. Torres-Martínez } \\ 
\normalsize  Department of Mathematics - UACM \\ 
\normalsize  \href{mailto:inocencio3@gmail.com}{inocencio3@gmail.com}
}
\date{} 
\renewcommand{\maketitlehookd}{%

\begin{abstract}

In this document we present the construction of a radial functions that have the objective of emulating the behavior of the radial basis function thin plate spline (TPS), which we will name as function TPS, we propose a way to partially and totally apply the fractional derivative to these functions to be used in interpolation problems, a proposal is presented to precondition the matrices generated in the interpolation problem using the $QR$ decomposition and finally   is proposed the form of a radial interpolant to be used when solving differential equations using the asymmetric collocation method.

\textbf{Key words:} Radial Basis Functions, Fractional Calculus, Fractional Derivative of Riemann-Liouville.
\end{abstract}
}

\begin{document}

% Print the title
\maketitle

%----------------------------------------------------------------------------------------
%	ARTICLE CONTENTS
%----------------------------------------------------------------------------------------

\section{Construction of Functions }

\subsection{Polynomials similar to the  function TPS}

The main idea in this section is to try to emulate the behavior of the radial basis function thin plate spline (TPS) , also known as polyharmonic spline \cite{gonzalez}:

\begin{eqnarray}\label{eq:01}
\small
\begin{array}{ll}
\Phi(r)=r^nlog(r), & n\in 2\nset{N},
\end{array}
\end{eqnarray}

in a  domain $ \Omega $ of the form

\begin{eqnarray*}
\small
\begin{array}{l}
\Omega=[0,1]\times[0,1],
\end{array}
\end{eqnarray*}

towards a domain of the form

\begin{eqnarray}\label{eq:02}
\small
\begin{array}{l}
\Omega_b=[0,b]\times[0,b],
\end{array}
\end{eqnarray}

to do this it must be taken into account that \eqref{eq:01}  satisfy

\begin{eqnarray*}
\small
\begin{array}{ll}
\Phi(0)=0, & \Phi'(0)=0,\\ 
\Phi(1)=0, & \Phi'(1)=1, 
\end{array}
\end{eqnarray*}

then for our purpose we look for a radial function $\Phi(r)$ such that

\begin{eqnarray}
\small
\begin{array}{ll}
\Phi(0)=0, & \Phi'(0)=0,
\end{array} \label{eq:03}\\
\small
\begin{array}{ll}
\Phi(b)=0, & \Phi'(b)=1,
\end{array} \label{eq:04}
\end{eqnarray}

to satisfy the conditions given in \eqref{eq:03} is taken a polynomial of the form

\begin{eqnarray*}
\small
\begin{array}{l}
\Phi(r)= a_1 r^{N+1}+ a_0r^N,
\end{array}
\end{eqnarray*}

where the coefficients $a_0$ and $a_1$ are determined by \eqref{eq:04}, and the value of $N$ will be given later, then

\begin{eqnarray*}
\small
\begin{array}{lll}
\Phi(b)=&a_1b^{N+1}+a_0b^N &=0,\\
\Phi'(b)=&a_1(N+1)b^N+a_0Nb^{N-1}&=1,
\end{array}
\end{eqnarray*}

in matrix form the previous system takes the form

\begin{eqnarray}
\small
\begin{array}{l}
\underbrace{
\begin{pmatrix}
b^{N+1} & b^N\\
(N+1)b^N & Nb^{N-1}
\end{pmatrix}}_{B}
\underbrace{
\begin{pmatrix}
a_1\\
a_0
\end{pmatrix}}_{a}= \underbrace{ \begin{pmatrix}
0\\1
\end{pmatrix}}_{c},
\end{array}\label{eq:05}
\end{eqnarray}

denoting by $det(B)$ to the determinant of the matrix $B$ from the previous system and doing a bit of algebra we obtain that

\begin{eqnarray*}
\small
\begin{array}{lll}
det(B)=-b^{2N} \neq 0 & \Leftrightarrow &b\neq 0,
\end{array}
\end{eqnarray*}

then the system \eqref{eq:05} It always has a solution, denoting now by $adj(B)$ the adjoint matrix of $B$ and using that

\begin{eqnarray}
\small
\begin{array}{l}
B^{-1}= \dfrac{1}{det(B)}adj(B),
\end{array} \label{eq:11}
\end{eqnarray}

we obtain that

\begin{eqnarray*}
\small
\begin{array}{l}
B^{-1}=\begin{pmatrix}
-Nb^{-1-N} &b^{-N}\\
(N+1)b^{-N} &-b^{1-N}
\end{pmatrix},
\end{array}
\end{eqnarray*}

and it is obtained as a solution to the system \eqref{eq:05}

\begin{eqnarray*}
\small
\begin{array}{l}
\begin{pmatrix}
a_1 \\ a_0
\end{pmatrix}=
\begin{pmatrix}
b^{-N} \\
-b^{1-N}
\end{pmatrix},
\end{array}
\end{eqnarray*}

with which we obtain the polynomial

\begin{eqnarray}
\small
\begin{array}{l}
\Phi(N,r)=b^{-N}r^{N+1}-b^{1-N}r^N,
\end{array}\label{eq:06}
\end{eqnarray}

by construction \eqref{eq:06} in the domain $\Omega_1$ fulfills that

\begin{eqnarray*}
\small
\begin{array}{l}
\Phi(N,r)=r^{N+1}-r^N \approx r^Nlog(r),
\end{array}
\end{eqnarray*}

\begin{figure}[!ht]
    \begin{tabular}{c}
    \begin{subfigure}{0.23\textwidth}
        \includegraphics[width=\textwidth]{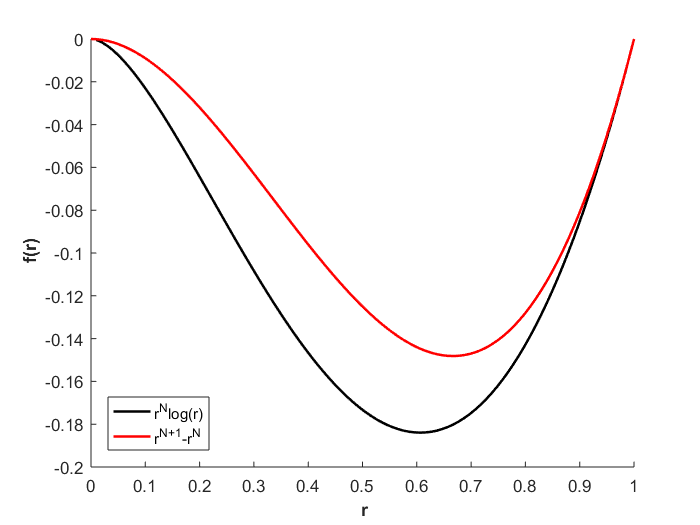}
        \caption{N=2}
    \end{subfigure}
    \begin{subfigure}[h]{0.23\textwidth}
        \includegraphics[width=\textwidth]{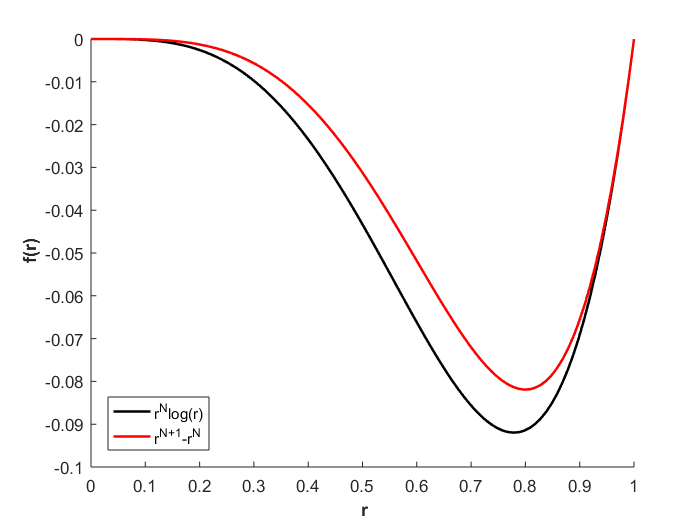}
        \caption{N=4}
    \end{subfigure}
    \\
    \begin{subfigure}{0.23\textwidth}
        \includegraphics[width=\textwidth]{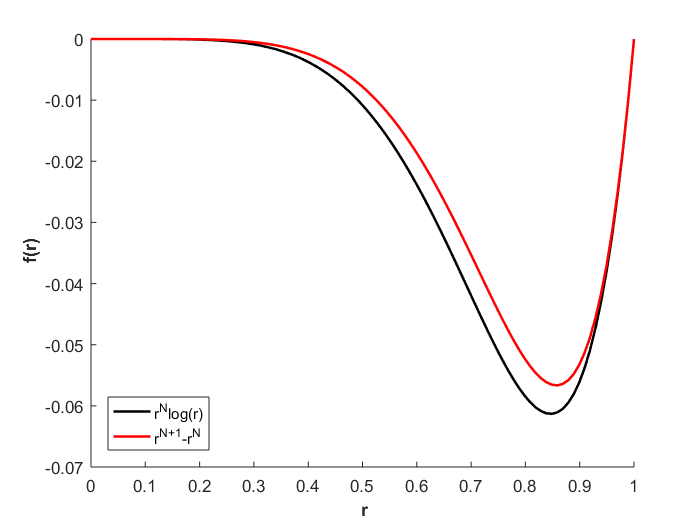}
        \caption{N=6}
    \end{subfigure}
    \end{tabular}
    \caption{They are presented with black and red the functions $r^Nlog(r)$ and $r^N-r^{N+1}$ respectively.}
\end{figure}

In the previous construction only two coefficients are used to perform the approximation of the function  TPS , to add one more coefficient we use the fact that \eqref{eq:01} in the domain $\Omega_1$ fulfills that

\begin{eqnarray*}
\small
\begin{array}{lll}
\Phi(0)=0, & \Phi'(0)=0, & \Phi''(0)=0\\ 
\Phi(1)=0, & \Phi'(1)=1, & \Phi''(1)=2n-1,
\end{array}
\end{eqnarray*}

then we look for a radial function $\Phi(r)$ such that

\begin{align}
&\small \begin{array}{lll}
\Phi(0)=0, & \Phi'(0)=0,&\Phi''(0)=0,
\end{array}\label{eq:08}\\
& \small \begin{array}{lll}
\Phi(b)=0, & \Phi'(b)=1,&\Phi''(b)=2N-1,
\end{array} \label{eq:09}
\end{align}

to satisfy \eqref{eq:08} we take the polynomial

\begin{eqnarray*}
\small
\begin{array}{l}
\Phi(r)=a_2r^{N+2}+a_1r^{N+1}+a_0r^N,
\end{array}
\end{eqnarray*}

on the other hand, to satisfy \eqref{eq:09} we arrived to the matrix system

\begin{eqnarray}
\scriptsize
\begin{array}{l}
\underbrace{
\begin{pmatrix}
b^{N+2}&b^{N+1} & b^N\\
(N+2)b^{N+1}&(N+1)b^N & Nb^{N-1}\\
(N+2)(N+1)b^N&(N+1)Nb^{N-1}&N(N-1)b^{N-2}
\end{pmatrix}}_{B}
\underbrace{
\begin{pmatrix}
a_2\\
a_1\\
a_0
\end{pmatrix}}_{a}= \underbrace{ \begin{pmatrix}
0\\1\\
2N-1
\end{pmatrix}}_{c}
,
\end{array}
\label{eq:10}
\end{eqnarray}

where to do a bit of algebra we get that

\begin{eqnarray*}
\small
\begin{array}{lll}
det(B)=-2b^{3N}\neq 0 & \Leftrightarrow &b\neq 0,
\end{array}
\end{eqnarray*}

and using \eqref{eq:11} we have to

\begin{eqnarray*}
\small
\begin{array}{l}
B^{-1}=\begin{pmatrix}
 \frac{1 }{2}N (N+1) b^{-2-N} & -Nb^{-1-N}  & \frac{1}{2} b^{-N}\\
 -N (N+2)b^{-1-N}  & (2N+1) b^{-N} & -b^{1-N} \\
 \frac{1}{2}\left(N^2+3N+2\right) b^{-N}  & -(N+1)b^{1-N}  & \frac{1}{2}b^{2-N}
\end{pmatrix},
\end{array}
\end{eqnarray*}

then the system \eqref{eq:10} has as a solution

\begin{eqnarray*}
\small
\begin{array}{l}
\begin{pmatrix}
a_2\\
a_1\\
a_0
\end{pmatrix}=
\begin{pmatrix}
 \frac{1}{2} (2 N-1)b^{-N} -N b^{-N-1} \\
 (2 N+1)b^{-N} -(2 N-1)b^{1-N}  \\
 \frac{1}{2}(2 N-1) b^{2-N} - (N+1) b^{1-N} 
\end{pmatrix},
\end{array}
\end{eqnarray*}

with which we get the polynomial

\begin{eqnarray} \label{eq:12}
\small
\begin{array}{ll}
\Phi(N,r)=&\left[\frac{1}{2} (2 N-1)b^{-N} -N b^{-N-1}\right] r^{N+2}  \\ 
&+\left[(2 N+1)b^{-N} -(2 N-1)b^{1-N} \right]r^{N+1}\\
&+\left[ \frac{1}{2}(2 N-1) b^{2-N} - (N+1) b^{1-N}  \right]r^N,
\end{array}
\end{eqnarray}

by construction \eqref{eq:12} in the domain $\Omega_1$ satisfy that

\begin{eqnarray*}
\small
\begin{array}{l}
\Phi(N,r)=-\frac{1}{2}r^{N+2}+2r^{N+1}-\frac{3}{2} r^N \approx r^Nlog(r),
\end{array}
\end{eqnarray*}

\newpage

\begin{figure}[!ht]
    \begin{tabular}{c}
    \begin{subfigure}{0.23\textwidth}
        \includegraphics[width=\textwidth]{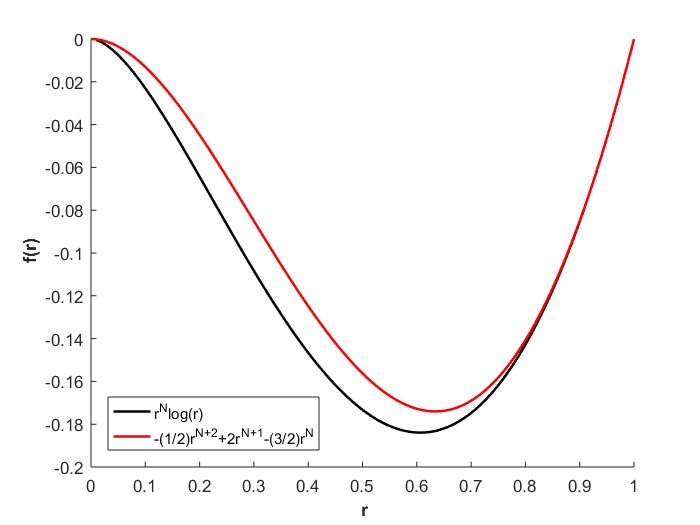}
        \caption{N=2}
    \end{subfigure}
    \begin{subfigure}[h]{0.23\textwidth}
        \includegraphics[width=\textwidth]{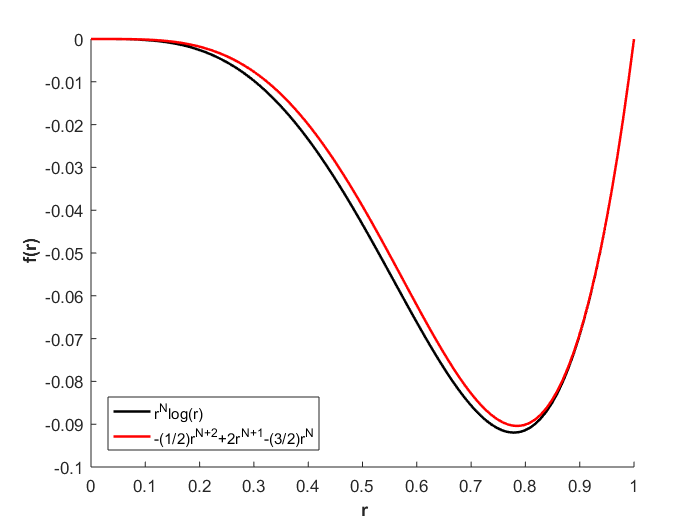}
        \caption{N=4}
    \end{subfigure}
    \\
    \begin{subfigure}{0.23\textwidth}
        \includegraphics[width=\textwidth]{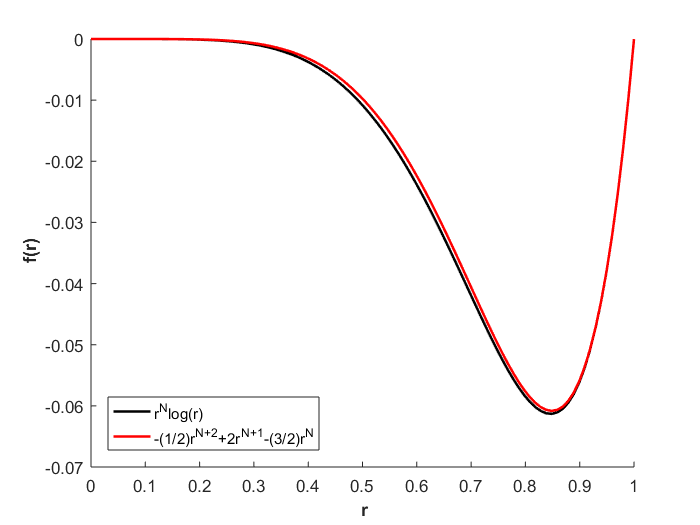}
        \caption{N=6}
    \end{subfigure}
    \end{tabular}
    \caption{They are presented with black and red the functions $r^Nlog(r)$ and $-\frac{1}{2}r^{N+2}+2r^{N+1}-\frac{3}{2} r^N$ respectively.}
\end{figure}

From the systems \eqref {eq:05} and \eqref{eq:10} it can be deduced that to construct a polynomial with $ n $ coefficients that approximates the function TPS , we need to consider the $ (n-1) $-derivatives from both the polynomial and the function TPS , but this path would make the coefficients $ a_i$'s become more  complicated in its expressions, to try to solve this problem in the next subsection we present an alternative that allows us to obtain an approximation to the function TPS where the coefficients $a_i's$ they are maintained in a simple way.

\subsubsection{Function false  TPS}

With the idea of later using the fractional derivative of polynomials \cite{oldham74} and keep retaining the behavior of the function TPS, which is zero at the extremes of the domain $\Omega_1$, we start with the idea of looking for a polynomial that becomes zero in the \com{extremes} with respect to the derivative, as the solution of the system \eqref{eq:05} with the above mentioned conditions it leads to the trivial solution., we used the polynomial involved in the system \eqref{eq:10} with a vector $ c $ of the form

\begin{eqnarray*}
\small
\begin{array}{l}
c=\begin{pmatrix}
0\\0\\
-c_0
\end{pmatrix},
\end{array}
\end{eqnarray*}

with $c_0>0$ and where the minus sign is included so that the solution has a convex behavior analogous to the function  TPS in the domain $\Omega_1$.

With the above the system \eqref{eq:10}  it can be rewritten as

\begin{eqnarray}
\scriptsize
\begin{array}{l}
\underbrace{
\begin{pmatrix}
b^{N+2}&b^{N+1} & b^N\\
(N+2)b^{N+1}&(N+1)b^N & Nb^{N-1}\\
(N+2)(N+1)b^N&(N+1)Nb^{N-1}&N(N-1)b^{N-2}
\end{pmatrix}}_{B}
\underbrace{
\begin{pmatrix}
a_2\\
a_1\\
a_0
\end{pmatrix}}_{a}= \underbrace{ \begin{pmatrix}
0\\0\\
-c_0
\end{pmatrix}}_{c},
\end{array}
\label{eq:14}
\end{eqnarray}

which has as solution

\begin{eqnarray*}
\small
\begin{array}{l}
\begin{pmatrix}
a_2\\
a_1\\
a_0
\end{pmatrix}=
\begin{pmatrix}
 -\frac{1}{2}c_0 b^{-N} \\
 c_0b^{1-N}  \\
 -\frac{1}{2}c_0b^{2-N} 
\end{pmatrix},
\end{array}
\end{eqnarray*}

with which we get the polynomial

\begin{eqnarray}\label{eq:15}
\small
\begin{array}{l}
\Phi(N,r)=-\dfrac{c_0}{2} b^{-N}r^{N+2}+c_0b^{1-N} r^{N+1}-\dfrac{c_0}{2}b^{2-N} r^N,
\end{array}
\end{eqnarray}

although in principle $ c_0 $ can be arbitrary,  later we propose a way to select it so that the coefficients of polynomial \eqref{eq:15} are maintain in a simple way, for the particular case $ c_0 = 4 $ we get the polynomial

\begin{eqnarray}\label{eq:16}
\small
\begin{array}{l}
\Phi(N,r)=-2 b^{-N}r^{N+2}+4b^{1-N} r^{N+1}-2b^{2-N} r^N,
\end{array}
\end{eqnarray}

the election of $ c_0 $ and the construction of \eqref{eq:16} guarantees that in the domain $\Omega_1$ fulfills that

\begin{eqnarray*}
\small
\begin{array}{l}
\Phi(N,r)=-2 r^{N+2}+4 r^{N+1}-2 r^N \approx r^Nlog(r).
\end{array}
\end{eqnarray*}

\begin{figure}[!ht]
    \begin{tabular}{c}
    \begin{subfigure}{0.23\textwidth}
        \includegraphics[width=\textwidth]{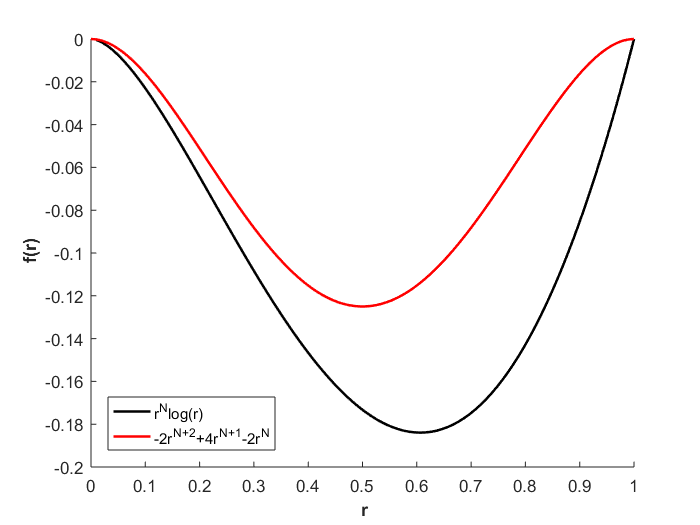}
        \caption{N=2}
    \end{subfigure}
    \begin{subfigure}[h]{0.23\textwidth}
        \includegraphics[width=\textwidth]{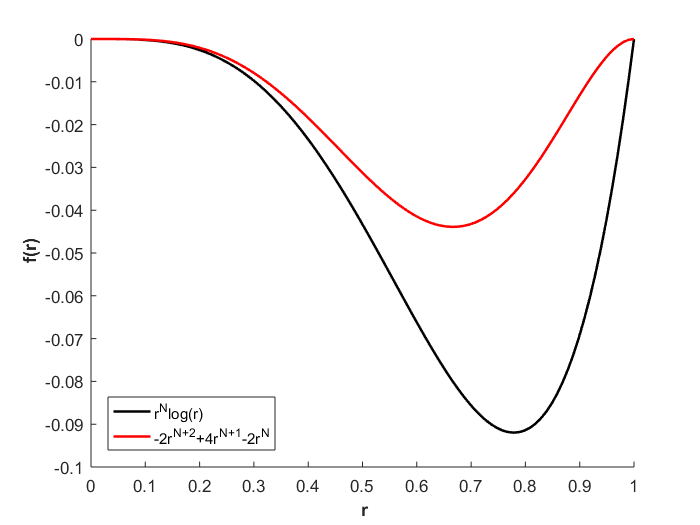}
        \caption{N=4}
    \end{subfigure}
    \\
    \begin{subfigure}{0.23\textwidth}
        \includegraphics[width=\textwidth]{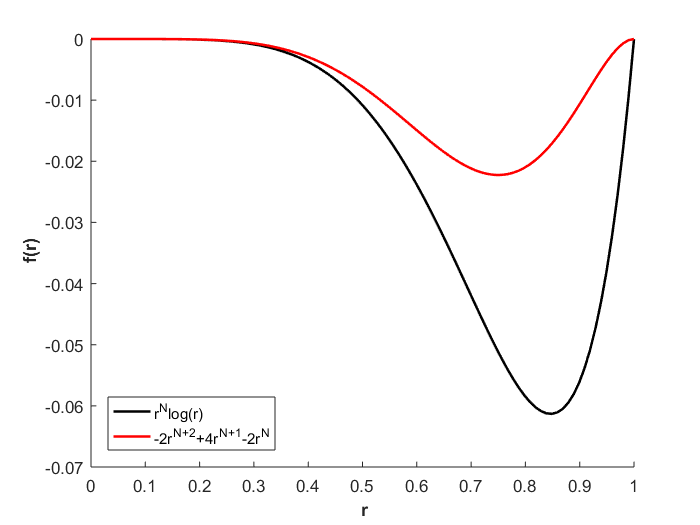}
        \caption{N=6}
    \end{subfigure}
    \end{tabular}
    \caption{They are presented with black and red the functions $r^Nlog(r)$ and $-2r^{N+2}+4r^{N+1}-2r^N$ respectively.}
\end{figure}

To improve the approximation we taken a small perturbation  $-\alpha$, with $\alpha\in[0,1)$, in the exponent of the term of greater power associated with a negative coefficient, modifying in turn the exponent of said coefficient with a value $+\alpha$, then we can define the function

\begin{eqnarray}\label{eq:17}
\footnotesize
\begin{array}{l}
\Phi(\alpha, N,r)=-2 b^{-N+\alpha}r^{N-\alpha+2}+4b^{1-N} r^{N+1}-2b^{2-N} r^N,
\end{array}
\end{eqnarray}

which in the domain $\Omega_1$ fulfills that
\begin{eqnarray}\label{eq:18}
\small
\begin{array}{l}
\Phi(\alpha,N,r)=-2r^{N-\alpha+2}+4r^{N+1}-2 r^N \approx r^Nlog(r).
\end{array}
\end{eqnarray}

The equation \eqref{eq:18} receives the name of false TPS while the equation \eqref{eq:17} receives the name of false  TPS generalized.

\begin{figure}[!ht]
    \begin{tabular}{c}
    \begin{subfigure}{0.23\textwidth}
        \includegraphics[width=\textwidth]{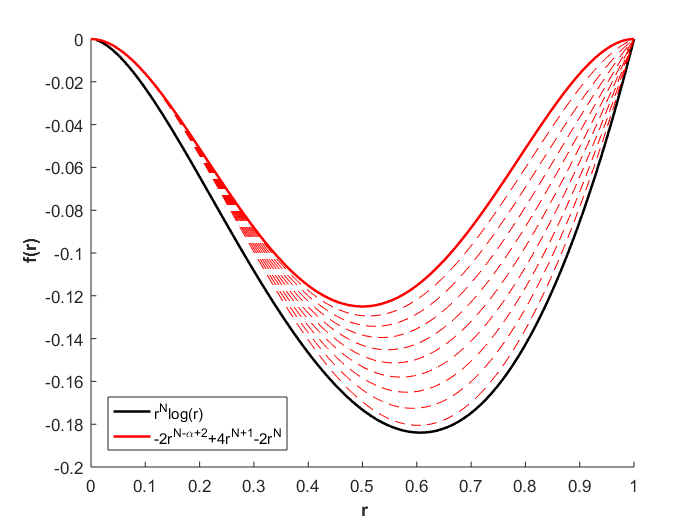}
        \caption{N=2}
    \end{subfigure}
    \begin{subfigure}[h]{0.23\textwidth}
        \includegraphics[width=\textwidth]{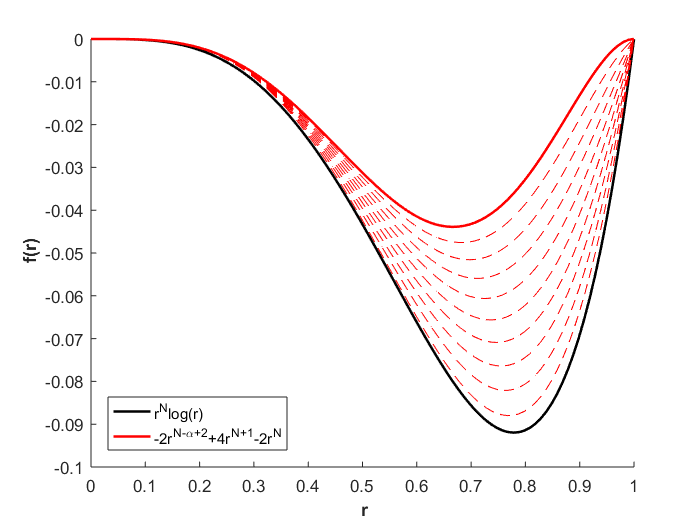}
        \caption{N=4}
    \end{subfigure}
    \\
    \begin{subfigure}{0.23\textwidth}
        \includegraphics[width=\textwidth]{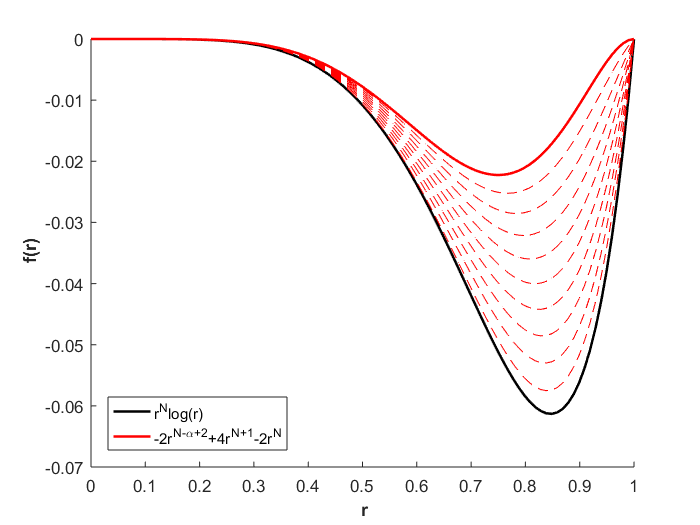}
        \caption{N=6}
    \end{subfigure}
    \end{tabular}
    \caption{They are presented with black and red the functions $r^Nlog(r)$ and $-2r^{N-\alpha+2}+4r^{N+1}-2r^N$, using different values of $\alpha$,  respectively. }
\end{figure}

\subsubsection{Generalizing the previous construction}

To generalize the idea used in the construction of the function false TPS, it seeks to generate a polynomial of the form

\begin{eqnarray*}
\small
\begin{array}{l}
\Phi(r)=a_3r^{N+3}+a_2r^{N+2}+a_1r^{N+1}+a_0r^N,
\end{array}
\end{eqnarray*}

that satisfies the conditions

\begin{eqnarray*}
\small
\begin{array}{llll}
\Phi(b)=0, & \Phi'(b)=0, &\Phi''(b)=0,& \Phi'''(b)=c_0,
\end{array}
\end{eqnarray*}

this generates a matrix system of the form.

\begin{eqnarray*}
\small
\begin{array}{l}
Ba=c,
\end{array}
\end{eqnarray*}

where

\begin{eqnarray*}
\small
\begin{array}{l}
c=\begin{pmatrix}
0\\0\\0\\c_0
\end{pmatrix},
\end{array}
\end{eqnarray*}

and

\begin{eqnarray*}
\small
\begin{array}{lll}
det(B)=12b^{4N}\neq 0 &\Leftrightarrow & b\neq 0,
\end{array}
\end{eqnarray*}

on the other hand using \eqref{eq:11} we have to

\begin{eqnarray*}
\small
\begin{array}{l}
B^{-1}=\begin{pmatrix}
B_3^{-1}&B_2^{-1} &B_1^{-1}&B_0^{-1}
\end{pmatrix},
\end{array}
\end{eqnarray*}

where $\set{B_i^{-1}}_{i=0}^3$ are the column vectors of the inverse matrix of $B$, with

\begin{eqnarray*}
\small
\begin{array}{l}
B_0^{-1}=\begin{pmatrix}
\frac{1}{6}b^{-N}\\
-\frac{1}{2}b^{1-N}\\
\frac{1}{2}b^{2-N}\\
-\frac{1}{6}b^{3-N}
\end{pmatrix},
\end{array}
\end{eqnarray*}

then the matrix system has like solution

\begin{eqnarray*}
\small
\begin{array}{l}
a=\begin{pmatrix}
\frac{c_0}{6}b^{-N}\\
-\frac{c_0}{2}b^{1-N}\\
\frac{c_0}{2}b^{2-N}\\
-\frac{c_0}{6}b^{3-N}
\end{pmatrix},
\end{array}
\end{eqnarray*}

with which we get the polynomial

\begin{eqnarray}
\small
\begin{array}{ll}
\Phi(N,r)=&\dfrac{c_0}{6}b^{-N}r^{N+3}-\dfrac{c_0}{2}b^{1-N}r^{N+2}\\
&+\dfrac{c_0}{2}b^{2-N} r^{N+1}-\dfrac{c_0}{6}b^{3-N}r^N,
\end{array}\label{eq:20}
\end{eqnarray}

denoting by $M=mcm(2,6)$, the lowest common multiple of the denominators present in the coefficients of \eqref{eq:20}, it defines

\begin{eqnarray*}
\small
\begin{array}{ll}
c_0=pM, & p\in \nset{Z}\setminus \set{0},
\end{array}
\end{eqnarray*}

where

\begin{eqnarray*}
\small
\begin{array}{l}
p=\left\{
\begin{array}{ll}
>0 , &\mbox{if $\Phi(2,r)$ is convex in $\Omega_1$} \\
<0,&\mbox{if $\Phi(2,r)$  is concave in $\Omega_1$}
\end{array}\right.,
\end{array}
\end{eqnarray*}

\newpage

to \eqref{eq:20} we can take  $c_0=18$, obtaining the polynomial

\begin{eqnarray}\label{eq:21}
\scriptsize
\begin{array}{l}
\Phi(N,r)=3b^{-N}r^{N+3}-9b^{1-N}r^{N+2}
+9b^{2-N} r^{N+1}-3b^{3-N}r^N,
\end{array}
\end{eqnarray}

due to the choice of $c_0$ and the way it is built \eqref{eq:21} we have to in the domain $\Omega_1$ fulfills that

\begin{eqnarray*}
\small
\begin{array}{l}
\Phi( N,r)=3r^{N+3}-9r^{N+2}+9 r^{N+1}-3r^N
  \approx r^Nlog(r),
\end{array}
\end{eqnarray*}

\begin{figure}[!ht]
    \begin{tabular}{c}
    \begin{subfigure}{0.23\textwidth}
        \includegraphics[width=\textwidth]{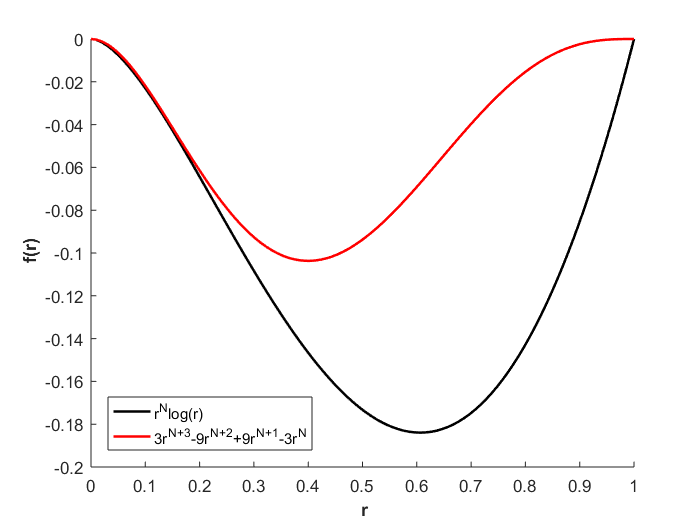}
        \caption{N=2}
    \end{subfigure}
    \begin{subfigure}[h]{0.23\textwidth}
        \includegraphics[width=\textwidth]{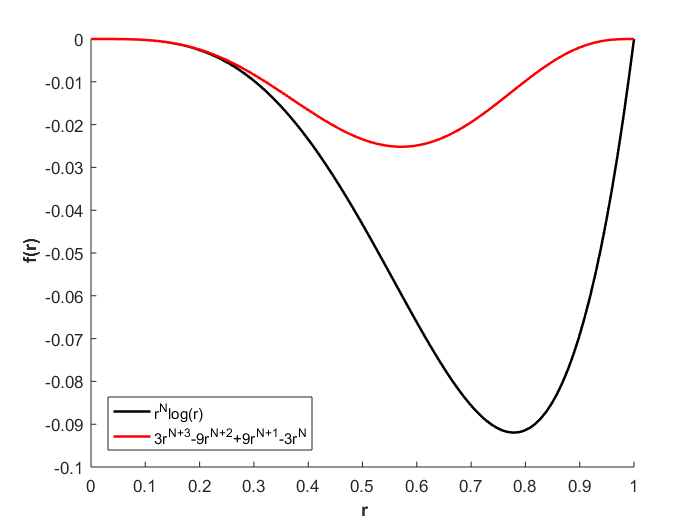}
        \caption{N=4}
    \end{subfigure}
    \\
    \begin{subfigure}{0.23\textwidth}
        \includegraphics[width=\textwidth]{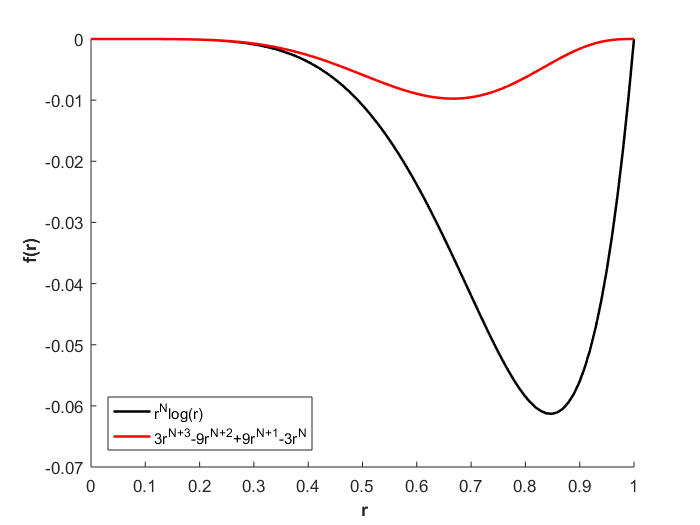}
        \caption{N=6}
    \end{subfigure}
    \end{tabular}
    \caption{They are presented with black and red the functions $r^Nlog(r)$ and $3r^{N+3}-9r^{N-\alpha+2}+9 r^{N+1}-3r^N$ respectively.}
\end{figure}

To improve the approximation take a small perturbation $-\alpha$, with $\alpha\in[0,1)$, in the exponent of the term of greater power associated with a negative coefficient, modifying in turn the exponent of said coefficient with a value $+\alpha$, then we can define the function

\begin{eqnarray}\label{eq:22}
\small
\begin{array}{ll}
\Phi(\alpha, N,r)=&3b^{-N}r^{N+3}-9b^{1-N+\alpha}r^{N-\alpha+2}\\
&+9b^{2-N} r^{N+1}-3b^{3-N}r^N,
\end{array}
\end{eqnarray}

which in the domain $\Omega_1$ fulfills that

\begin{eqnarray}\label{eq:23}
\scriptsize
\begin{array}{l}
\Phi(\alpha, N,r)=3r^{N+3}-9r^{N-\alpha+2}+9 r^{N+1}-3r^N
 \approx r^Nlog(r),
\end{array}
\end{eqnarray}

\newpage

\begin{figure}[!ht]
    \begin{tabular}{c}
    \begin{subfigure}{0.23\textwidth}
        \includegraphics[width=\textwidth]{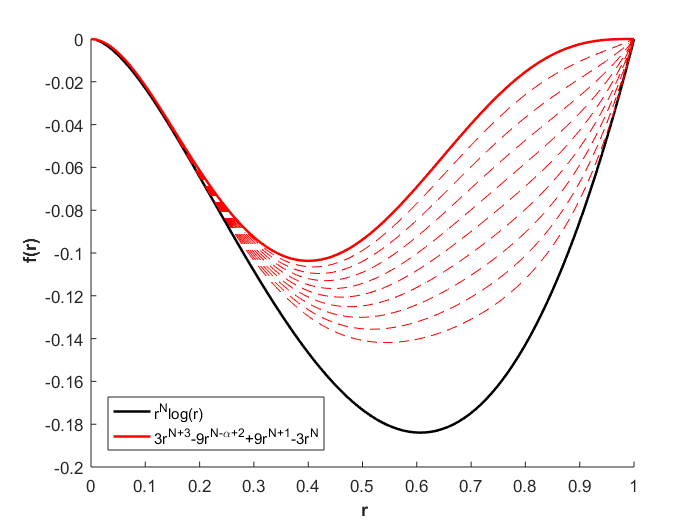}
        \caption{N=2}
    \end{subfigure}
    \begin{subfigure}[h]{0.23\textwidth}
        \includegraphics[width=\textwidth]{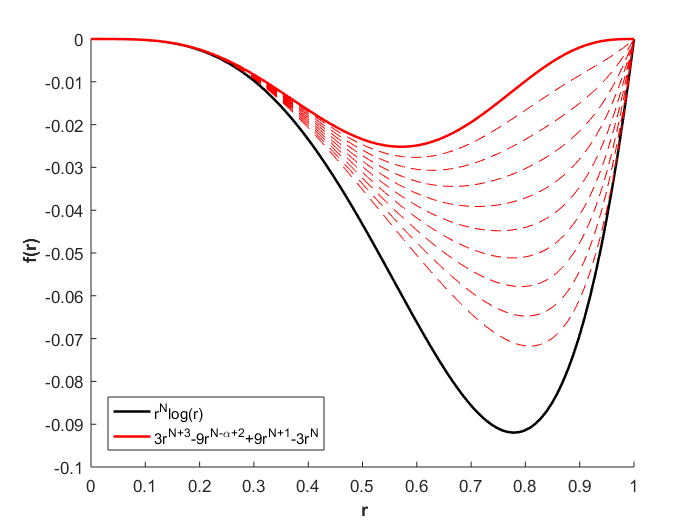}
        \caption{N=4}
    \end{subfigure}
    \\
    \begin{subfigure}{0.23\textwidth}
        \includegraphics[width=\textwidth]{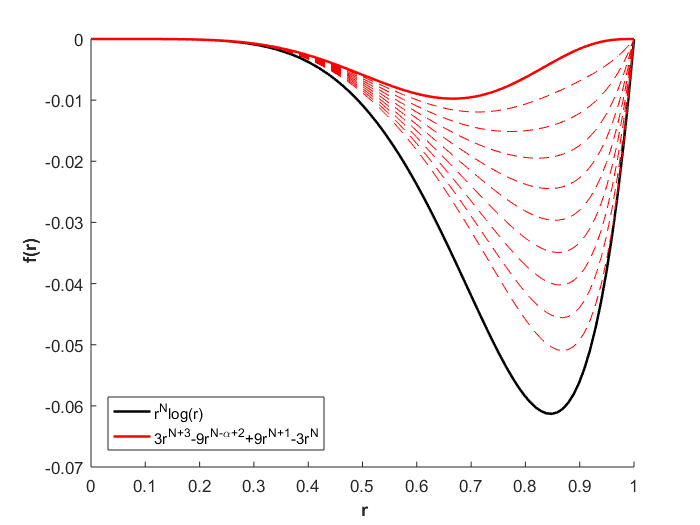}
        \caption{N=6}
    \end{subfigure}
    \end{tabular}
    \caption{They are presented with black and red the functions $r^Nlog(r)$ and $3r^{N+3}-9r^{N-\alpha+2}+9 r^{N+1}-3r^N$, using different values of $\alpha$,  respectively.}
\end{figure}

Of the way in which is constructed the polynomial  \eqref{eq:06} it can be generalized by changing the vector $c$  by the vector

\begin{eqnarray*}
\small
\begin{array}{l}
c'=\begin{pmatrix}
0\\c_0
\end{pmatrix},
\end{array}
\end{eqnarray*}

getting

\begin{eqnarray}\label{eq:25}
\small
\begin{array}{l}
\Phi(N,r)=c_0b^{-N}r^{N+1}-c_0b^{1-N}r^N,
\end{array}
\end{eqnarray}

taking the particular case $c_0=1$ and to improve the approximation take a small perturbation  $-\alpha$, with $\alpha\in[0,1)$, in the exponent of the term of greater power associated with a negative coefficient, modifying in turn the exponent of said coefficient with a value $+\alpha$, then we can define the function

\begin{eqnarray}\label{eq:07}
\small
\begin{array}{l}
\Phi(\alpha,N,r)=b^{-N}r^{N+1}-b^{1-N+\alpha}r^{N-\alpha},
\end{array}
\end{eqnarray}

which in the domain $ \ Omega_1 $ fulfills that

\begin{eqnarray}\label{eq:39}
\small
\begin{array}{l}
\Phi(\alpha,N,r)=r^{N+1}-r^{N-\alpha}\approx r^Nlog(r),
\end{array}
\end{eqnarray}

\newpage

\begin{figure}[!ht]
    \begin{tabular}{c}
    \begin{subfigure}{0.23\textwidth}
        \includegraphics[width=\textwidth]{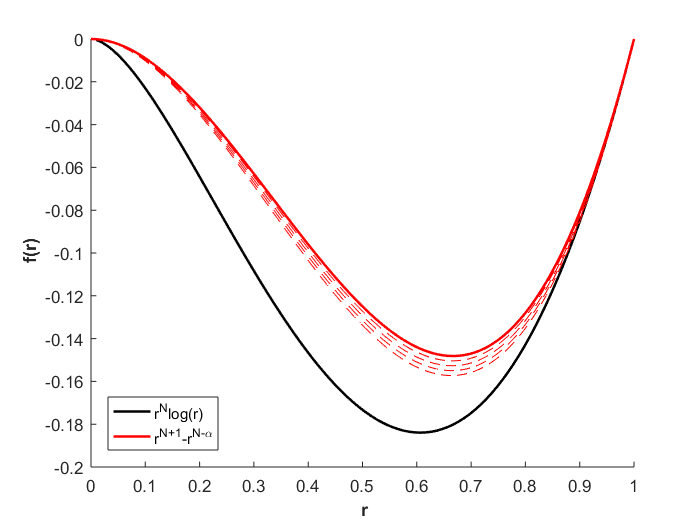}
        \caption{N=2}
    \end{subfigure}
    \begin{subfigure}[h]{0.23\textwidth}
        \includegraphics[width=\textwidth]{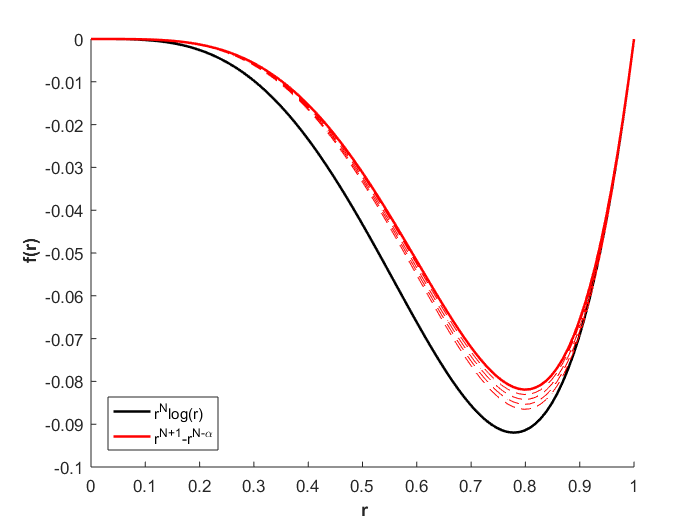}
        \caption{N=4}
    \end{subfigure}
    \\
    \begin{subfigure}{0.23\textwidth}
        \includegraphics[width=\textwidth]{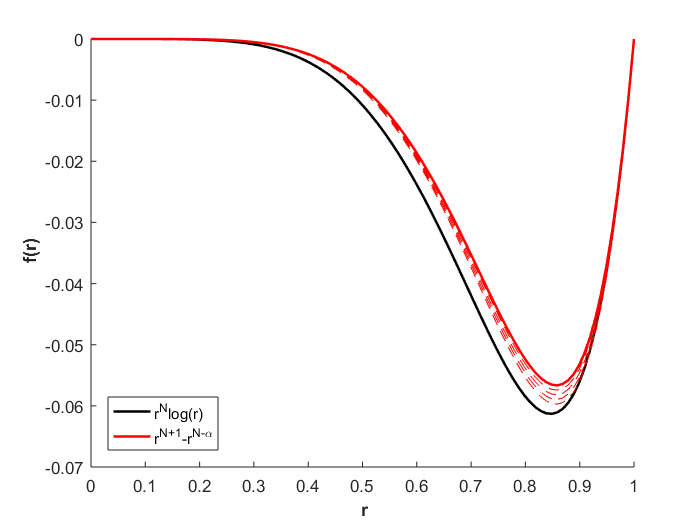}
        \caption{N=6}
    \end{subfigure}
    \end{tabular}
    \caption{They are presented with black and red the functions $r^Nlog(r)$ and $r^{N+1}-r^{N-\alpha}$,  using different values of $\alpha$,  respectively.}
\end{figure}

\subsection{Radial functions similar to the function TPS }

The functions \eqref{eq:07}, \eqref{eq:17} and \eqref{eq:22} behave similarly to the TPS function in the domain $\Omega_1$, but our purpose is to obtain radial functions \cite{gonzalez,wendland} that satisfied the previously mentioned, to solve this we impose the restrictions

\begin{eqnarray}\label{eq:13}
\small
\begin{array}{lll}
N\notin \nset{N} & y & \left( N-\alpha \right) \notin \nset{N},
\end{array}
\end{eqnarray}

donde $N>0$ y $\alpha \in [0,1)$. 

From now on we will take that all the functions used will have implicitly the restrictions given in \eqref{eq:13} unless otherwise mentioned.

Imposing the restrictions \eqref{eq:13} to the polynomials \eqref{eq:07}, \eqref{eq:17} and \eqref{eq:22} it is guaranteed that we have radial functions that behave similarly to the function TPS , to visualize this we  choose  the false function TPS and allow \eqref{eq:01} take rational values obtaining the following graphs

\newpage

\begin{figure}[!ht]
    \begin{tabular}{c}
    \begin{subfigure}{0.23\textwidth}
        \includegraphics[width=\textwidth]{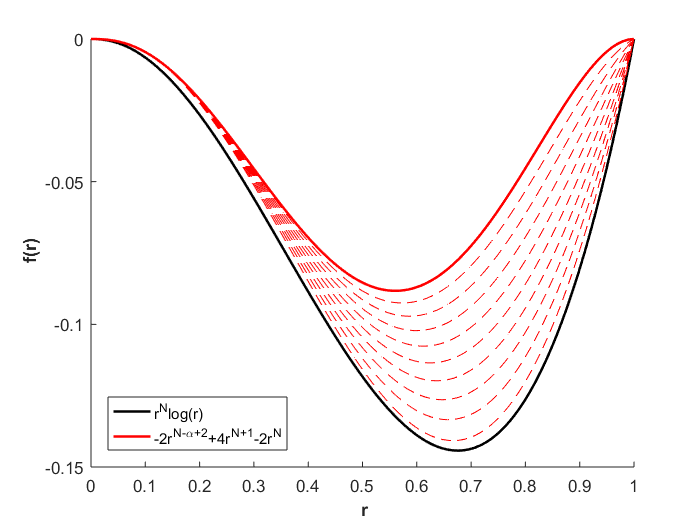}
        \caption{N=2.55}
    \end{subfigure}
    \begin{subfigure}[h]{0.23\textwidth}
        \includegraphics[width=\textwidth]{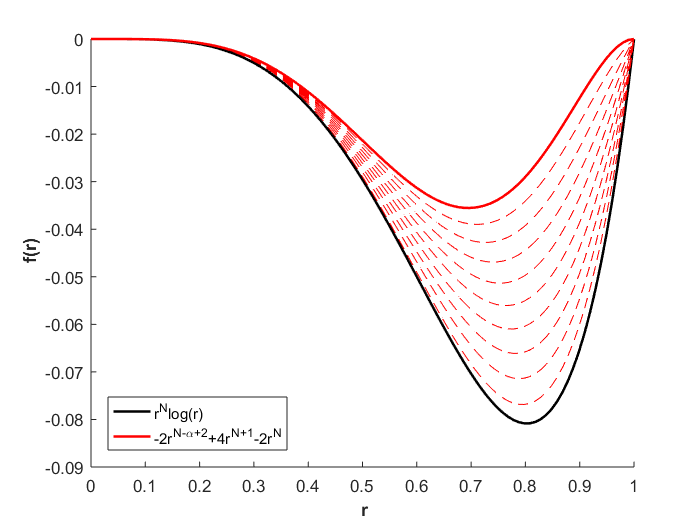}
        \caption{N=4.55}
    \end{subfigure}
    \\
    \begin{subfigure}{0.23\textwidth}
        \includegraphics[width=\textwidth]{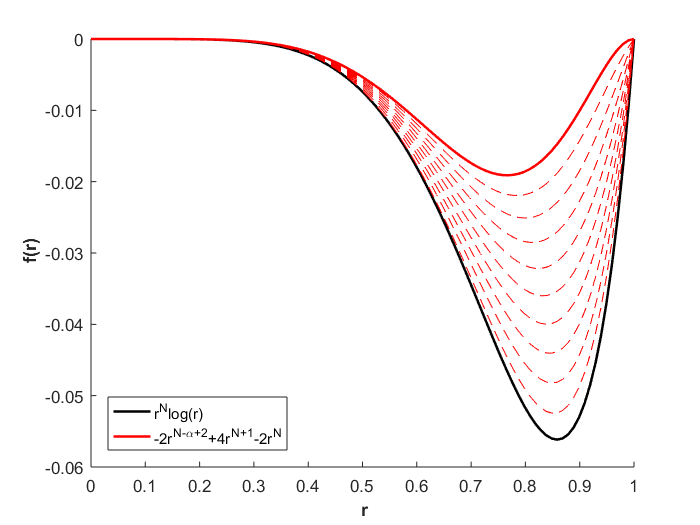}
        \caption{N=6.55}
    \end{subfigure}
    \end{tabular}
    \caption{They are presented with black and red the functions $r^Nlog(r)$ and $-2r^{N-\alpha+2}+4r^{N+1}-2r^N$, using different values of $\alpha$, respectively.}
\end{figure}

\subsubsection{Conditionally positive definite functions}
We start the next section giving a definition and a theorem \cite{wendland}  that will be very useful later

\begin{definition}\label{def:01}
A function $\phi : [0,\infty) \to \nset{R}$ which belongs to $C[0,\infty) \cap C ^\infty (0,\infty)$ and it satisfies

\begin{eqnarray}
\small
\begin{array}{ll}
(-1)^l\phi^{(l)}(r)\geq 0 , & \forall l\in \nset{N},
\end{array}
\end{eqnarray}

where $r>0$, it is called a completely monotone functions in $[0,\infty)$.

\end{definition}

\begin{theorem}
(\textbf{Michelli}) Suppose that $\phi \in C [0 ,\infty ) \cap C^\infty (0 ,\infty)$ is given. Then the function $\Phi = \phi ( \norm{\ \cdot \ }^2 )$ is
radial and conditionally positive defined order $m$ in $\nset{R}^d$ for all $d$ if and only if $( - 1)^m\phi^{(m)}$ is completely monotone in $[0,\infty)$.
\end{theorem}

We now consider the following example

\begin{example}
Suppose that $\phi$ is given by

\begin{eqnarray*}
\small
\begin{array}{ll}
\phi(r)=(-1)^{\lceil \beta/2 \rceil}r^{\beta/2}, & 0<\beta \notin \nset{N},
\end{array}
\end{eqnarray*}

where $r>0$, then

\begin{eqnarray*}
\small
\begin{array}{ll}
\phi^{(1)}(r)=&(-1)^{\lceil \beta/2 \rceil}\dfrac{\beta}{2}r^{(\beta/2)-1},\\
\phi^{(2)}(r)=&(-1)^{\lceil \beta/2 \rceil}\dfrac{\beta}{2}\left(\dfrac{\beta}{2}-1 \right) r^{(\beta/2)-2}, \\
\vdots& \\
\phi^{(l)}(r)=&(-1)^{\lceil \beta/2 \rceil}\dfrac{\beta}{2}\left(\dfrac{\beta}{2}-1 \right)\cdots\left(\dfrac{\beta}{2}-l+1 \right) r^{(\beta/2)-l},
\end{array}
\end{eqnarray*}

rewriting the last expression

\begin{eqnarray*}
\small
\begin{array}{l}
\ds \phi^{(l)}(r)= (-1)^{\lceil \beta/2 \rceil} \left[\prod_{k=1}^l \left(\dfrac{\beta}{2}-k+1 \right) \right]r^{(\beta/2)-l},
\end{array}
\end{eqnarray*}

then

\begin{eqnarray*}
\footnotesize{
\begin{array}{l}
\ds (-1)^{\lceil \beta/2 \rceil}\phi^{(\lceil \beta/2 \rceil)}(r)= \left[\prod_{k=1}^{\lceil \beta/2 \rceil} \left(\dfrac{\beta}{2}-k+1 \right) \right]r^{(\beta/2)-\lceil \beta/2 \rceil},
\end{array}
}
\end{eqnarray*}

which implies that

\begin{eqnarray*}
\small
\begin{array}{lll}
(-1)^{\lceil \beta/2 \rceil}\phi^{(\lceil \beta/2 \rceil)} \geq 0 & \Leftrightarrow & -1\leq \dfrac{\beta}{2}-k , \ \forall k \\
& \Leftrightarrow & -1\leq \dfrac{\beta}{2}- \left\lceil \dfrac{\beta}{2} \right\rceil \\
& \Leftrightarrow & \left\lceil \dfrac{\beta}{2} \right\rceil \leq \dfrac{\beta}{2}+1,
\end{array}
\end{eqnarray*}

with which $(-1)^{\lceil \beta/2 \rceil}\phi^{(\lceil \beta/2 \rceil)}$ is completely monotone, it should be noted that $m=\lceil \beta/2 \rceil$ is the smallest number for which this is fulfilled. Since $\beta$ is not a natural number, $\phi$ it is not a polynomial, and therefore the powers

\begin{eqnarray*}
\small
\begin{array}{ll}
\Phi(x)=(-1)^{\lceil \beta/2 \rceil}\norm{x}^{\beta}, & 0<\beta \notin \nset{N},
\end{array}
\end{eqnarray*}

they are strictly conditionally positive definite of order $\lceil \beta/2 \rceil$ and radials in $\nset{R}^d$ for all $d$.

\end{example}

A conditionally positive definite function of order $m$ , is also conditionally positive definite function of order $l \geq m $. It is also true that if a function is conditionally
positive definite of order $m$ in $\nset{R}^d $, then it is conditionally positive definite of order $m$ in $\nset{R}^k$ , for $k \leq d$ \cite{gonzalez}.

With the previous example we have the false function TPS

\begin{eqnarray*}
\small
\begin{array}{l}
\Phi(\alpha,N,r)=-2r^{N-\alpha+2}+4r^{N+1}-2 r^N, 
\end{array}
\end{eqnarray*}

is conditionally positive definite of order 

\begin{eqnarray*}
\small
\begin{array}{l}
\left\lceil \dfrac{N-\alpha+2}{2} \right\rceil
\end{array}.
\end{eqnarray*}

\section{Interpolation with Radial Functions}

A function $\Phi: \nset{R}^d \to \nset{R}$ is called radial, if there is a function $\phi : [0 ,\infty ) \to \nset{R}$ such that

\begin{eqnarray}\label{eq:19}
\small
\begin{array}{l}
\Phi(x)=\phi \left(\norm{x}\right),
\end{array}
\end{eqnarray}

where $\norm{ \ \cdot \ }$ is the Euclidean norm in $\nset{R}^d $.

Given a set of values $\set{( x_j ,u_j )}_{j=1}^{N_p} $, where $( x_j ,u_j )\in \Omega \times \nset{R}$  with $\Omega \subseteq \nset{R}^d$, an interpolant is a function $\sigma : \Omega \to \nset{R}$ such that

\begin{eqnarray}\label{eq:24}
\small
\begin{array}{ll}
\sigma(x_j ) = u_j ,& j \in \set{1,...,N_p}.
\end{array}
\end{eqnarray}

When is used a radial function  $\Phi$ conditionally positive definite, an interpolant of the form is proposed \cite{gonzalez,wendland}

\begin{eqnarray}\label{eq:27}
\small
\begin{array}{l}
\ds \sigma(x)=\sum_{j=1}^{N_p}\lambda_j \Phi(x-x_j) + \sum_{k=1}^Q \beta_k p_k(x),
\end{array}
\end{eqnarray}

where $Q = dim\left( \nset{P}_{m-1} \left(\nset{R}^d \right)\right)$ and $\set{p_k}_{k=1}^Q$ it is a base for $\nset{P}_{m-1} \left(\nset{R}^d \right)$. The Interpolation conditions \eqref{eq:24} are completed with the  moment conditions

\begin{eqnarray}\label{eq:28}
\small
\begin{array}{ll}
\ds{\sum_{j=1}^{N_p} \lambda_jp_k(x_j)=0} , & k\in \set{1,\cdots,Q}.
\end{array}
\end{eqnarray}

Solve the problem of interpolation \eqref{eq:24} using the interpolant \eqref{eq:27} together with the moment conditions  \eqref{eq:28} is equivalent to solving the linear system

\begin{eqnarray}\label{eq:29}
\small
\begin{array}{l}
\underbrace{
\begin{pmatrix}
A&P\\
P^T&0
\end{pmatrix}}_{G}
\underbrace{
\begin{pmatrix}
\lambda \\
\beta
\end{pmatrix}}_{\Lambda}=  \underbrace{\begin{pmatrix}
u\\0
\end{pmatrix}}_{U},
\end{array}
\end{eqnarray}

where $A$ and $P$ are matrices of $(N_p \times N_p)$ y $(N_p \times Q)$ respectively, whose components are

\begin{eqnarray}\label{eq:30}
\small
\begin{array}{ll}
A_{jk}=\Phi(x_j-x_k), & P_{jk}=p_k(x_j).
\end{array}
\end{eqnarray}

The one that a function $\Phi$ be conditionally positive definite of order $m$ , it can be interpreted as the matrix $A$ of components $A_{jk} = \Phi( x_j - x_k )$ is positive definite in the space of
vectors $c\in \nset{R}^{N_p}$ such that

\begin{eqnarray}\label{eq:26}
\small
\begin{array}{ll}
\ds{ \sum_{j=1}^{N_p}c_j p_k (x_j ) = 0}, & k\in \set{1,\cdots, Q}.
\end{array}
\end{eqnarray}

In this sense, $ A $ is positive definite in the vector space
 $c$ that are \com{perpendiculars} to
the polynomials. So, if in \eqref{eq:27} the function $\Phi$ is conditionally positive definite
of order $m$ and the set of centers $\set{x_j}_{j=1}^{N_p}$ contains a unisolvent subset, then the interpolation problem will have a solution (the condition of unisolvency is to ensure the uniqueness) \cite{gonzalez}.

\subsection{Examples with Radial Functions}

Defining a domain

\begin{eqnarray*}
\small
\begin{array}{l}
\Omega_{a,b}:=[a,b]\times[a,b],
\end{array}
\end{eqnarray*}

and using the function

\begin{eqnarray}\label{eq:31}
\small
\begin{array}{l}
u(x,y)=\dfrac{ \sin(8(x+y))+ \cos(8(x-y))+4}{35},
\end{array}
\end{eqnarray}

with the following distribution of Halton type inner nodes over the domain $\Omega_{0.28,1.48}$

\begin{figure}[!ht]
\includegraphics[width=0.45\textwidth]{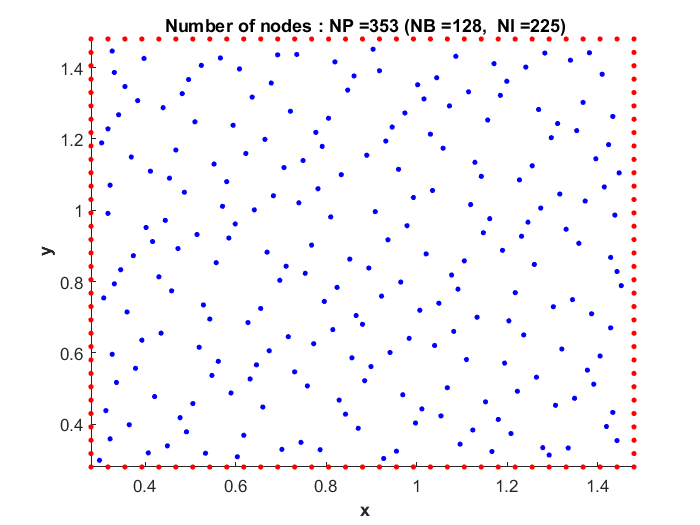}
\caption{Nodes used for the interpolation problem, where $NB$ and $NI$ are the boundary and interior nodes respectively.
}\label{fig:01}
\end{figure}

we get that the graph of \eqref{eq:31} is given by

\begin{figure}[!ht]
\includegraphics[width=0.45\textwidth]{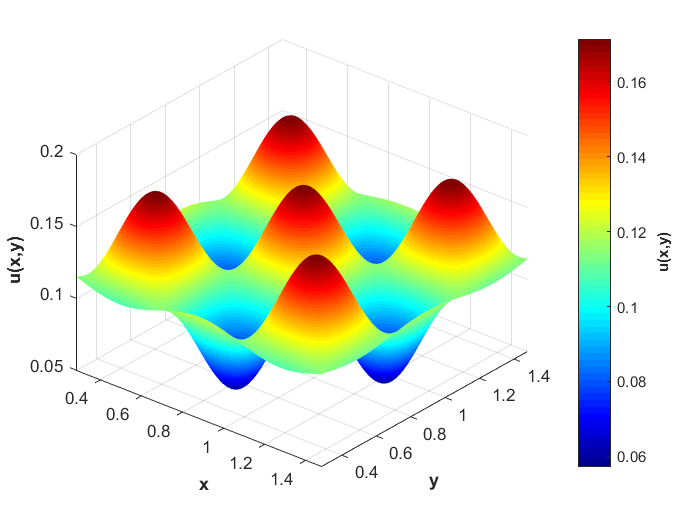}
\caption{Graph of the equation \eqref{eq:31}.}
\end{figure}

Then to carry out the interpolation problem a set of values is generated $\set{u(x_i,y_i)}_{i=1}^{N_p}$ and it is taken  $\alpha \in [0,1)$, here the option of using $ \alpha $ for fixed values is presented,  although it can also be used by looking for a value that minimizes the error.

Denoting by $\sigma_i=\sigma(x_i,y_i)$ and $u_i=u(x_i,y_i)$, then the error that we will use will be the root of the mean square error given by

\begin{eqnarray}\label{eq:43}
RMSE:=\sqrt{ \dfrac{1}{N_p}\sum_{i=1}^{N_p}\left(u_i-\sigma_i \right)^2},
\end{eqnarray}

denoting by $cond(G)$ the condition number of the  matrix $ G $ given in \eqref{eq:29}, the following examples are presented

\begin{itemize}
\item[•] Using the  false function TPS generalized

\begin{eqnarray*}
\small
\begin{array}{l}
\Phi(\alpha, N,r)=-2 b^{-N+\alpha}r^{N-\alpha+2}+4b^{1-N} r^{N+1}-2b^{2-N} r^N,
\end{array}
\end{eqnarray*}

taking $N=3.22$, then to use the interpolant \eqref{eq:27} it defines

\begin{eqnarray*}
\small
\begin{array}{l}
m=\max \set{\left\lceil \dfrac{N-0+2}{2} \right\rceil,\left\lceil \dfrac{N-1+2}{2} \right\rceil }=3,
\end{array}
\end{eqnarray*}

obtaining the following results

\begin{eqnarray*}
\footnotesize
\begin{array}{ccc}
\alpha& RMSE  & cond(G) \\ \hline
0.0	&	1.0914474968265677e-11	&	1.2378900238537703e7	\\
 0.1	&	3.5419628199123710e-11	&	1.2735545436520265e7	\\
 0.2	&	1.8453287132781697e-11	&	1.3129133985928234e7	\\
 0.3	&	1.6149600468160241e-11	&	1.3562189431763934e7	\\
 0.4	&	7.3342258102292580e-11	&	1.4037225997380065e7	\\
 0.5	&	1.0865604380572651e-11	&	1.4556544783656418e7	\\
 0.6	&	7.9374223725322362e-11	&	1.5121867662453054e7	\\
 0.7	&	4.2596214388606749e-11	&	1.5733707496013103e7	\\
 0.8	&	1.7399138544255595e-10	&	1.6390320782187937e7	\\
 0.9	&	7.7515629759964534e-11	&	1.7086025186181001e7		
\end{array}
\end{eqnarray*}

\item[•]Using the radial function

\begin{eqnarray*}
\small
\begin{array}{ll}
\Phi(\alpha, N,r)=&3b^{-N}r^{N+3}-9b^{1-N+\alpha}r^{N-\alpha+2}\\
&+9b^{2-N} r^{N+1}-3b^{3-N}r^N,
\end{array}
\end{eqnarray*}

taking $N=2.55$, then to use the interpolant \eqref{eq:27} it defines

\begin{eqnarray*}
\small
\begin{array}{l}
m=\left\lceil \dfrac{N+3}{2} \right\rceil=3,
\end{array}
\end{eqnarray*}

obtaining the following results

\begin{eqnarray*}
\footnotesize
\begin{array}{ccc}
\alpha& RMSE  & cond(G) \\ \hline
0.0	&	4.7298423751436141e-12	&	7.7628112197229778e5	\\
 0.1	&	6.0250788793601656e-11	&	9.9014734372436150e5	\\
 0.2	&	1.2464988456346324e-11	&	1.2568932888810737e6	\\
 0.3	&	4.6859085196101791e-11	&	3.1146550821565916e6	\\
 0.4	&	1.0027959159749939e-11	&	1.9080520595207722e6	\\
 0.5	&	2.0267447654438357e-11	&	2.2522535371709852e6	\\
 0.6	&	1.4284160002323938e-11	&	2.6389644941461636e6	\\
 0.7	&	5.2356342185294705e-12	&	3.0161926367063778e6	\\
 0.8	&	7.3474961415774206e-12	&	3.3928655142875575e6	\\
 0.9	&	1.8963726382781025e-11	&	3.7524027369971084e6		
\end{array}
\end{eqnarray*}

\end{itemize}

\section{Fractional Derivative}

The perturbations $-\alpha$ previously used have a structure similar to the fractional derivative of Riemann-Liouville \cite{oldham74,miller93}, which in its unified form with the fractional integral of Riemann-Liouville
  \cite{hilfer00} is given by

\begin{eqnarray}
\scriptsize
\begin{array}{l}
\ifr{a}{}{D}{x}{\alpha}f(x)= \left\{
\begin{array}{cl}
\displaystyle{ \frac{1}{\gam{-\alpha}}\int_a^x(x-t)^{-\alpha -1}f(t)dt,} &\re{\alpha}<0\\ \\
\displaystyle{ \frac{1}{\gam{n-\alpha}}\dfrac{d^n}{dx^n}\int_a^x(x-t)^{n-\alpha -1}f(t)dt,} & \re{\alpha}> 0 
\end{array}\right.
\end{array}
, \nonumber \\
\label{eq:32}
\end{eqnarray}

where $n=\lfloor\re{\alpha}\rfloor+1$. 

For a monomial given by $f(x)=x^s$, the fractional derivative of Riemann-Liouville takes the form

\begin{eqnarray}\label{eq:33}
\small
\begin{array}{l}
\ifr{0}{}{D}{x}{\alpha}x^s= \dfrac{\gam{s+1}}{\gam{s-\alpha+1}}x^{s-\alpha},
\end{array}
\end{eqnarray}

then to implement the fractional derivative to the radial functions \eqref{eq:17} and \eqref{eq:22} is taken

\begin{eqnarray*}
\small
\begin{array}{l}
r^{N-\alpha+2} \to \ifr{0}{}{D}{r}{\alpha}r^{N+2},
\end{array}
\end{eqnarray*}

getting the functions

\begin{eqnarray}\label{eq:34}
\small{
\begin{array}{ll}
\Phi(\alpha, N,r)=&-2 b^{-N+\alpha}\ifr{0}{}{D}{r}{\alpha}r^{N+2}+4b^{1-N} r^{N+1}\\
&-2b^{2-N}r^N,
\end{array}
}
\end{eqnarray}

\begin{eqnarray}\label{eq:40}
\small
\begin{array}{ll}
\Phi(\alpha, N,r)=&3b^{-N}r^{N+3}-9b^{1-N+\alpha}\ifr{0}{}{D}{r}{\alpha}r^{N+2}\\
&+9b^{2-N} r^{N+1}-3b^{3-N}r^N.
\end{array}
\end{eqnarray}

\subsection{Preconditioning of system}

Before continuing we must note that in the previous examples the condition number of the matrices obtained is too high
, also that the linear system \eqref{eq:29} generated to carry out the interpolation can be written compactly as

\begin{eqnarray*}
\small
\begin{array}{l}
G\Lambda =U,
\end{array}
\end{eqnarray*}

where $G$ is a matrix of $(N_p+Q)\times(N_p+Q)$, similarly $\Lambda$ and $U$ they are column vectors of $(N_p+Q)$ entries, then to try to solve the problem of having a condition number too high we propose to use the decomposition $QR$ \cite{plato2003concise} of the matrix $G$

\begin{eqnarray*}
\small
\begin{array}{l}
G=QR,
\end{array}
\end{eqnarray*}

and change the linear system \eqref{eq:29} by the equivalent linear system

\begin{eqnarray}\label{eq:40}
\small
\begin{array}{l}
G_M \Lambda=\left[ \left( HR\right)^{-1}G \right]\Lambda =(HR)^{-1}U,
\end{array}
\end{eqnarray}

where

\begin{eqnarray*}
\small
\begin{array}{ll}
H_{ij}= Q_{ij}+\dfrac{1}{2^n},&n\in \nset{N},
\end{array}
\end{eqnarray*}

taking the value of $ n $ in such a way that it is satisfied

\begin{eqnarray*}
\small
\begin{array}{ll}
cond\left(G_M \right)\leq M, & M<cond(G).
\end{array}
\end{eqnarray*}

In the following examples will be used  the linear system  \eqref{eq:40} using $M=10$.

\subsection{Examples with Fractional Derivative implemented partially}

Using again the equation  \eqref{eq:31}, the distribution of nodes of the Figure \ref{fig:01} and the set of values $\set{u_i}_{i=1}^{N_p}$, also how it is used the definition of fractional derivative given in \eqref{eq:32} is taken  $\alpha \in (-1,1)$, here is presented the option to use $ \ alpha $ for fixed values although it can also be used looking for a value that minimizes the error. The following examples are presented

\begin{itemize}
\item[•] Using the radial function

\begin{eqnarray*}
\small{
\begin{array}{ll}
\Phi(\alpha, N,r)=&-2 b^{-N+\alpha}\ifr{0}{}{D}{r}{\alpha}r^{N+2}+4b^{1-N} r^{N+1}\\
&-2b^{2-N}r^N,
\end{array}
}
\end{eqnarray*}

taking $N=3.22$, then to use the interpolant \eqref{eq:27} it defines

\begin{eqnarray*}
\small
\begin{array}{l}
m=\max \set{\left\lceil \dfrac{N-(-1)+2}{2} \right\rceil,\left\lceil \dfrac{N-1+2}{2} \right\rceil }=4,
\end{array}
\end{eqnarray*}

obtaining the following results

\begin{eqnarray*}
\footnotesize
\begin{array}{ccc}
\alpha& RMSE & cond(G_M) \\ \hline
 -0.9	&	7.0843389855493818e-9	&	6.8962078566406095	\\
 -0.8	&	7.9627100410608106e-9	&	6.0569588639096121	\\
 -0.7	&	6.9023958352949166e-9	&	5.1963557874926538	\\
 -0.6	&	3.2035903493421426e-9	&	4.3434965781152250	\\
 -0.5	&	3.4972838502603867e-9	&	8.3052228894332512	\\
 -0.4	&	1.0582732272271552e-9	&	6.0419636627818498	\\
 -0.3	&	2.1157007230952478e-9	&	4.2210575043174874	\\
 -0.2	&	2.1709843703594850e-9	&	6.9200465802306219	\\
 -0.1	&	1.8997640101087966e-10	&	8.5984545751152321	\\
  0.0	&	4.4378079985347997e-11	&	4.5437369066354600	\\
  0.1	&	5.6566327154318037e-10	&	8.9586468651466049	\\
  0.2	&	1.0593792560351835e-9	&	3.8743124718533268	\\
  0.3	&	1.0254581655063807e-8	&	4.5313967438105802	\\
  0.4	&	2.2137081038571393e-8	&	4.4140777166867364	\\
  0.5	&	2.5667950365750804e-8	&	3.8965071350163751	\\
  0.6	&	4.9843375967404835e-8	&	3.9787909923665166	\\
  0.7	&	5.7211155585308675e-8	&	7.1988316914494419	\\
  0.8	&	1.6752839287329331e-7	&	4.9090050755603283	\\
  0.9	&	1.9199368201557929e-7	&	8.8940930876438280	
\end{array}
\end{eqnarray*}

\item[•] Using the radial function

\begin{eqnarray*}
\small
\begin{array}{ll}
\Phi(\alpha, N,r)=&3b^{-N}r^{N+3}-9b^{1-N+\alpha}\ifr{0}{}{D}{r}{\alpha}r^{N+2}\\
&+9b^{2-N} r^{N+1}-3b^{3-N}r^N.
\end{array}
\end{eqnarray*}

taking $N=2.55$, then to use the interpolant \eqref{eq:27} it defines

\begin{eqnarray*}
\small
\begin{array}{l}
m=\left\lceil \dfrac{N+3}{2} \right\rceil=3,
\end{array}
\end{eqnarray*}

obtaining the following results

\begin{eqnarray*}
\footnotesize
\begin{array}{ccc}
\alpha& RMSE & cond(G_M) \\ \hline
 -0.9	&	1.5636401280910224e-9	&	8.8176620286344551	\\
 -0.8	&	5.8873758675912041e-9	&	6.1360976947294148	\\
 -0.7	&	5.4092063701652622e-8	&	5.0433564177708901	\\
 -0.6	&	4.3374345369719424e-9	&	6.3196876248218699	\\
 -0.5	&	8.1067055207989931e-9	&	6.8841681638364154	\\
 -0.4	&	3.3878992828354837e-9	&	9.7405683180427651	\\
 -0.3	&	8.0989820535485463e-7	&	4.4361589923520315	\\
 -0.2	&	1.3506585465600225e-9	&	4.0611147290345038	\\
 -0.1	&	4.3215604269889866e-10	&	4.1752341663354136	\\
  0.0	&	1.3702328833012031e-12	&	4.3432720452267706	\\
  0.1	&	1.4806716960271588e-7	&	6.7663396247438836	\\
  0.2	&	6.3294009384505283e-9	&	8.4599217905754376	\\
  0.3	&	1.3531823751896382e-8	&	4.2422896467663094	\\
  0.4	&	2.7243124773336195e-8	&	9.2532652527731329	\\
  0.5	&	9.0038145734016248e-9	&	7.3628765046197051	\\
  0.6	&	1.1028135893299789e-8	&	5.9925266668660280	\\
  0.7	&	4.7634941304921479e-8	&	5.7933459310608697	\\
  0.8	&	5.0410074306971711e-8	&	8.7060490429955451	\\
  0.9	&	2.5122663838290705e-8	&	4.9469813263343037			
\end{array}
\end{eqnarray*}

\end{itemize}

\subsection{Examples with Fractional Derivative}

Because the previous examples where partial fractional derivative is implemented did not present any problem to carry out the problem of interpolation, is proceeded to implement the fractional derivative in its entirety. To implement the fractional derivative to the radial functions \eqref{eq:17} and \eqref{eq:22}  is taken

\begin{eqnarray*}
\small
\begin{array}{ll}
b^{s} r^{t}  \to  b^{s+\alpha}\ifr{0}{}{D}{r}{\alpha}r^{t},
\end{array}
\end{eqnarray*}

getting the functions

\begin{eqnarray}\label{eq:35}
\small{
\begin{array}{ll}
\Phi(\alpha, N,r)=&-2 b^{-N+\alpha}\ifr{0}{}{D}{r}{\alpha}r^{N+2}+4b^{1-N+\alpha}\ifr{0}{}{D}{r}{\alpha} r^{N+1}\\
&-2b^{2-N+\alpha}\ifr{0}{}{D}{r}{\alpha}r^N,
\end{array}
}
\end{eqnarray}

\begin{eqnarray}
\small
\begin{array}{ll}
\Phi(\alpha, N,r)=&3b^{-N+\alpha}\ifr{0}{}{D}{r}{\alpha}r^{N+3}-9b^{1-N+\alpha}\ifr{0}{}{D}{r}{\alpha}r^{N+2}\\
&+9b^{2-N+\alpha}\ifr{0}{}{D}{r}{\alpha} r^{N+1}-3b^{3-N+\alpha}\ifr{0}{}{D}{r}{\alpha}r^N.
\end{array}
\end{eqnarray}

Using  again the equation \eqref{eq:31}, the distribution of nodes of the Figure \ref{fig:01} and the set of values $\set{u_i}_{i=1}^{N_p}$, also how it is used the definition of fractional derivative given in \eqref{eq:32} is taken  $\alpha \in (-1,1)$, here is presented the option to use $ \ alpha $ for fixed values although it can also be used looking for a value that minimizes the error. The following examples are presented

\begin{itemize}
\item[•]Using the radial function

\begin{eqnarray*}
\small{
\begin{array}{ll}
\Phi(\alpha, N,r)=&-2 b^{-N+\alpha}\ifr{0}{}{D}{r}{\alpha}r^{N+2}+4b^{1-N+\alpha}\ifr{0}{}{D}{r}{\alpha} r^{N+1}\\
&-2b^{2-N+\alpha}\ifr{0}{}{D}{r}{\alpha}r^N,
\end{array}
}
\end{eqnarray*}

taking $N=3.22$, then to use the interpolant \eqref{eq:27} it defines

\begin{eqnarray*}
\small
\begin{array}{l}
m=\max \set{\left\lceil \dfrac{N-(-1)+2}{2} \right\rceil,\left\lceil \dfrac{N-1+2}{2} \right\rceil }=4,
\end{array}
\end{eqnarray*}

obtaining the following results

\begin{eqnarray*}
\footnotesize
\begin{array}{ccc}
\alpha& RMSE & cond(G_M)\\ \hline
 -0.9	&	5.8908349572860871e-5	&	5.5737193415777879	\\
 -0.8	&	1.8509839045727311e-7	&	9.3325125607453536	\\
 -0.7	&	1.9925355752856633e-8	&	6.5758364796013336	\\
 -0.6	&	1.1431508773757680e-8	&	9.9240318018805915	\\
 -0.5	&	1.0448728429386133e-9	&	9.3638113635796358	\\
 -0.4	&	9.6214768828011250e-10	&	6.1275981202468977	\\
 -0.3	&	1.8528894507976145e-10	&	5.7042585623510025	\\
 -0.2	&	9.0446894698987107e-11	&	3.8298337481043938	\\
 -0.1	&	3.3158975735329348e-11	&	4.3052843086581811	\\
  0.0	&	4.4378079985347997e-11	&	4.5437369066354600	\\
  0.1	&	3.6007027267234390e-12	&	6.3594701648028167	\\
  0.2	&	3.0493121013947708e-12	&	4.0907320406004128	\\
  0.3	&	3.7303122196736812e-12	&	5.5527905168744915	\\
  0.4	&	2.1017398075672669e-12	&	7.7639923195233767	\\
  0.5	&	2.8642668952531466e-12	&	4.4302484678083367	\\
  0.6	&	1.7151280653557004e-12	&	9.1424671423519843	\\
  0.7	&	3.2023309197061125e-10	&	6.2845169502316427	\\
  0.8	&	3.5711150057257600e-11	&	8.9905660498588915	\\
  0.9	&	5.9638735502978287e-11	&	8.2997802368672993		
\end{array}
\end{eqnarray*}

\item[•] Using the radial function

\begin{eqnarray*}
\small
\begin{array}{ll}
\Phi(\alpha, N,r)=&3b^{-N+\alpha}\ifr{0}{}{D}{r}{\alpha}r^{N+3}-9b^{1-N+\alpha}\ifr{0}{}{D}{r}{\alpha}r^{N+2}\\
&+9b^{2-N+\alpha}\ifr{0}{}{D}{r}{\alpha} r^{N+1}-3b^{3-N+\alpha}\ifr{0}{}{D}{r}{\alpha}r^N.
\end{array}
\end{eqnarray*}

taking $N=2.55$, then to use the interpolant \eqref{eq:27} it defines

\begin{eqnarray*}
\small
\begin{array}{l}
m=\left\lceil \dfrac{N+3}{2} \right\rceil=3,
\end{array}
\end{eqnarray*}

obtaining the following results

\begin{eqnarray*}
\footnotesize
\begin{array}{ccc}
\alpha& RMSE& cond(G_M)\\ \hline
 -0.9	&	1.3021967648637492e-10	&	8.3206817644804421	\\
 -0.8	&	1.3155695012794909e-10	&	4.8835303263483159	\\
 -0.7	&	2.3228970398179848e-11	&	9.7249628177357774	\\
 -0.6	&	2.2909145676738737e-11	&	4.9463326639225471	\\
 -0.5	&	1.1222370660857585e-11	&	9.8783558056367120	\\
 -0.4	&	4.2011829661131651e-12	&	3.8074930989810314	\\
 -0.3	&	1.6013645309492655e-12	&	9.4611712552245724	\\
 -0.2	&	3.8906902215986858e-12	&	7.7025873509331255	\\
 -0.1	&	5.6908853021616975e-12	&	4.9472748604993955	\\
  0.0	&	1.3702328833012031e-12	&	4.3432720452267706	\\
  0.1	&	3.0439857307602263e-12	&	4.8455204940138428	\\
  0.2	&	2.7638402095506380e-12	&	9.5434630857671525	\\
  0.3	&	7.0712507905215194e-12	&	3.9806740182868281	\\
  0.4	&	8.7277730326344969e-8	&	9.2824895482898651	\\
  0.5	&	1.6293963338135556e-11	&	4.9231085065761375	\\
  0.6	&	1.3205804213762410e-13	&	7.4458485925621094	\\
  0.7	&	7.0923741669988860e-14	&	8.8386103435077299	\\
  0.8	&	3.1225485178923595e-14	&	8.2851951713076737	\\
  0.9	&	9.7113169582163791e-15	&	8.3415405059857033	
\end{array}
\end{eqnarray*}

\end{itemize}

\subsubsection{A change in the interpolant}

In the previous sections we use the interpolator given by \eqref{eq:27} where $Q = dim\left( \nset{P}_{m-1} \left(\nset{R}^d \right)\right)$, this causes the value of $ Q $ to grow considerably, take for example a polynomial in $\nset{R}^2$ of degree $4$ which makes that $Q$ be equal to $15$, considering that sometimes \com{less is more} is changes the polynomial present in \eqref{eq:27} by a radial polynomial obtaining the following interpolant

\begin{eqnarray}\label{eq:36}
\small
\begin{array}{l}
\ds \sigma(x)=\sum_{j=1}^{N_p}\lambda_j\Phi(x-x_j)+\sum_{k=0}^Q\beta_kr^{k}(x),
\end{array}
\end{eqnarray}

where now $Q = dim\left( \nset{P}_{m-1} \left(\nset{R} \right)\right)$, with which the moment conditions takes the form

\begin{eqnarray}\label{eq:37}
\small
\begin{array}{ll}
\ds{ \sum_{j=1}^{N_p}\lambda_j r^{k}=0, }& k\in \set{0,1,\cdots,Q},
\end{array}
\end{eqnarray}

With which if we take now a radial polynomial in
 $\nset{R}$ of degree $4$ the value of $ Q $ would be equal to $5$. The following examples are presented with the interpolant mentioned above

\begin{itemize}
\item[•] Using the radial function

\begin{eqnarray*}
\small{
\begin{array}{ll}
\Phi(\alpha, N,r)=&-2 b^{-N+\alpha}\ifr{0}{}{D}{r}{\alpha}r^{N+2}+4b^{1-N+\alpha}\ifr{0}{}{D}{r}{\alpha} r^{N+1}\\
&-2b^{2-N+\alpha}\ifr{0}{}{D}{r}{\alpha}r^N,
\end{array}
}
\end{eqnarray*}

taking $N=3.22$, then to use the interpolant \eqref{eq:27} it defines

\begin{eqnarray*}
\small
\begin{array}{l}
m=\max \set{\left\lceil \dfrac{N-(-1)+2}{2} \right\rceil,\left\lceil \dfrac{N-1+2}{2} \right\rceil }=4,
\end{array}
\end{eqnarray*}

obtaining the following results

\begin{eqnarray*}
\footnotesize
\begin{array}{ccc}
\alpha& RMSE&cond(G_M)\\ \hline
 -0.9	&	1.5889551253924140e-5	&	9.4845006930392888	\\
 -0.8	&	6.3987623204658784e-7	&	6.9266877560988895	\\
 -0.7	&	5.5122415022730115e-8	&	4.5581991313406451	\\
 -0.6	&	4.9609699309717730e-9	&	3.9022712999605416	\\
 -0.5	&	2.8606758000057324e-9	&	9.7575009733381535	\\
 -0.4	&	8.5378019419285573e-10	&	4.1162090749498113	\\
 -0.3	&	1.6402729631796765e-10	&	4.9854689034480550	\\
 -0.2	&	3.2385258469302527e-10	&	7.1369658368082733	\\
 -0.1	&	3.0112296271981445e-11	&	8.0787947550106871	\\
  0.0	&	2.3959266689019608e-11	&	3.8057013377543933	\\
  0.1	&	7.9731842066765081e-12	&	4.2653333610861530	\\
  0.2	&	5.2174158637825498e-12	&	4.0587788299867702	\\
  0.3	&	4.7773108979483316e-12	&	4.4679461226648147	\\
  0.4	&	2.1669594332569641e-12	&	9.4197027627527756	\\
  0.5	&	3.0119577908331812e-12	&	7.8564160299507737	\\
  0.6	&	3.2593186988301638e-12	&	6.9239100545265257	\\
  0.7	&	9.3963218133542650e-11	&	5.5156116996767359	\\
  0.8	&	1.7046483976602697e-11	&	4.2959022223965544	\\
  0.9	&	5.3198946788454349e-11	&	9.8919840455960859			
\end{array}
\end{eqnarray*}

\item[•] Using the radial function

\begin{eqnarray*}
\small
\begin{array}{ll}
\Phi(\alpha, N,r)=&3b^{-N+\alpha}\ifr{0}{}{D}{r}{\alpha}r^{N+3}-9b^{1-N+\alpha}\ifr{0}{}{D}{r}{\alpha}r^{N+2}\\
&+9b^{2-N+\alpha}\ifr{0}{}{D}{r}{\alpha} r^{N+1}-3b^{3-N+\alpha}\ifr{0}{}{D}{r}{\alpha}r^N.
\end{array}
\end{eqnarray*}

taking $N=2.55$, then to use the interpolant \eqref{eq:27} it defines

\begin{eqnarray*}
\small
\begin{array}{l}
m=\left\lceil \dfrac{N+3}{2} \right\rceil=3,
\end{array}
\end{eqnarray*}

obtaining the following results

\begin{eqnarray*}
\footnotesize
\begin{array}{ccc}
\alpha& RMSE & cond(G_M) \\ \hline
 -0.9	&	1.0308678210497606e-10	&	6.8773854919881963	\\
 -0.8	&	7.2570454505231712e-11	&	5.5595347475522914	\\
 -0.7	&	1.5060722203198524e-11	&	7.3991587307494582	\\
 -0.6	&	2.4352127937860650e-11	&	4.0837782175908206	\\
 -0.5	&	2.4156954952439390e-11	&	7.4970533998625086	\\
 -0.4	&	4.4217156255200803e-12	&	3.8051504395575053	\\
 -0.3	&	8.3505968376411197e-12	&	5.4227744487381342	\\
 -0.2	&	4.7016062967983923e-12	&	7.6611909313834730	\\
 -0.1	&	1.8858726776411597e-12	&	7.2630611733591328	\\
  0.0	&	2.4993866322564443e-12	&	4.0482647647287502	\\
  0.1	&	5.9517756220989523e-11	&	4.7794223621267342	\\
  0.2	&	1.0009489420367950e-11	&	8.4226146559807606	\\
  0.3	&	9.4604822587979618e-12	&	7.5397704723868832	\\
  0.4	&	1.0034655564359237e-9	&	4.4732515594629589	\\
  0.5	&	1.9499567422389998e-11	&	4.6521910301042295	\\
  0.6	&	1.5720325480500878e-13	&	8.1449698846275034	\\
  0.7	&	6.4128676128549415e-14	&	4.3687895324727934	\\
  0.8	&	2.6723781539337429e-14	&	4.9392428716518344	\\
  0.9	&	2.7181313200709048e-14	&	5.9990758313730330		
\end{array}
\end{eqnarray*}

\end{itemize}

\section{Asymmetrical collocation}

The interpolation technique presented above can also be applied to the solution of differential equations
 \cite{gonzalez}. Assuming we have a domain $\Omega\subset \nset{R}^d$ and the problem

\begin{eqnarray}\label{eq:44}
\small
\left\{
\begin{array}{ll}
\mathcal{L}u = f,& \Omega,\\
\mathcal{B}u = g, &\partial \Omega,
\end{array}\right.
\end{eqnarray}

where $f$ and $g$ are functions given ,  $\mathcal{L}$ and $\mathcal{B}$ linear differential operators, and $u$ is the solution to find.

Before continuing we will make a change in the interpolant \eqref{eq:36} that will help us avoid discontinuities due to the application of the operators $\mathcal{L}$ and $\mathcal{B}$. Denoting by $ord(\mathcal{L})$ the order of the differential operator $\mathcal{L}$, we define

\begin{eqnarray*}
\small
\begin{array}{c}
q=\max\set{ord(\mathcal{L}),ord(\mathcal{B})},
\end{array}
\end{eqnarray*}

defining now

\begin{eqnarray*}
o=\left\{
\begin{array}{cl}
q-1,& \mbox{si }q>0\\
0,& \mbox{si }q\leq 0
\end{array}\right.,
\end{eqnarray*}

then the interpolant \eqref{eq:36} it can be rewritten as

\begin{eqnarray}\label{eq:45}
\small
\begin{array}{l}
\ds \sigma(x)=\sum_{j=1}^{N_p}\lambda_j\Phi(x-x_j)+\beta_0 +\sum_{k=1}^Q\beta_kr^{k+o}(x),
\end{array}
\end{eqnarray}

where $Q = dim\left( \nset{P}_{m-1} \left(\nset{R} \right)\right)$, then the moment  conditions take the form

\begin{eqnarray}\label{eq:46}
\small
\begin{array}{c}
\ds{ \sum_{j=1}^{N_p}\lambda_j p_1=\sum_{j=1}^{N_p}\lambda_j =0},\\
\begin{array}{ll}
\ds{ \sum_{j=1}^{N_p}\lambda_j p_{k+1}=\sum_{j=1}^{N_p}\lambda_j r^{k+o}=0, }& k\in \set{1,\cdots,Q},
\end{array}
\end{array}
\end{eqnarray}

finally to the restrictions given in \eqref{eq:13} we must add a more restriction given by

\begin{eqnarray}\label{eq:47}
\small
\left\{
\begin{array}{ll}
N >q+\alpha, & \mbox{si }q>0,\\
N>\alpha, & \mbox{si }q\leq 0.
\end{array}\right.
\end{eqnarray}

When replacing the interpolant \eqref{eq:45} in the system \eqref{eq:44} is obtained

\begin{eqnarray*}
\small
\left\{
\begin{array}{ll}
\mathcal{L}\sigma = f,& \Omega,\\
\mathcal{B}\sigma = g, &\partial \Omega,
\end{array}\right.
\end{eqnarray*}

with which the linear system is obtained

\begin{eqnarray}
\small
\begin{array}{l}
\underbrace{
\begin{pmatrix}
\mathcal{L} A& \mathcal{L} P\\
\mathcal{B} A& \mathcal{B} P\\
P^T&0
\end{pmatrix}}_{G}
\underbrace{
\begin{pmatrix}
\lambda \\
\\
\beta
\end{pmatrix}}_{\Lambda}=  \underbrace{\begin{pmatrix}
f\\g\\0
\end{pmatrix}}_{U},
\end{array}
\end{eqnarray}

where $\mathcal{L} A$, $\mathcal{B}A$, $\mathcal{L} P$, $\mathcal{B}P$ and $P$   are matrices of $(N_I \times N_p)$, $((N_p-N_I) \times N_p)$, $(N_I \times Q)$, $((N_p-N_I) \times Q)$ and $(N_p \times Q)$ respectively, whose components are

\begin{eqnarray*}
\small
\begin{array}{c}
\begin{array}{ll}
\mathcal{L} A_{jk}=\mathcal{L}\Phi(x_j-x_k), &\mathcal{L} P_{jk}=\mathcal{L}p_k(x_j),\\
\mathcal{B}A_{jk}=\mathcal{B}\Phi(x_j-x_k),&\mathcal{B} P_{jk}=\mathcal{B}p_k(x_j),
\end{array}\\
\begin{array}{l}
 P_{jk}=p_k(x_j),
\end{array}
\end{array}
\end{eqnarray*}

\subsection{Examples with Fractional Derivative}

The form of the interpolant \eqref{eq:45} it will be very useful to solve differential equations in radial form. Taking the definition of fractional derivative of Caputo \cite{oldham74}

\begin{eqnarray}\label{eq:41}
\scriptsize
\begin{array}{ll}
\ifr{a}{C}{D}{x}{\alpha}f(x)=\ds{ \frac{1}{\gam{n-\alpha}}\int_a^x(x-t)^{n-\alpha -1}\dfrac{d^n}{dt^n}f(t)dt}, & \re{\alpha}> 0 ,
\end{array}
\end{eqnarray}

we can build the next differential operator

\begin{eqnarray}\label{eq:42}
\small
\begin{array}{l}
\mathcal{L}=\left(\ifr{0}{C}{D}{r}{2+\beta} + \dfrac{1}{r} \ifr{0}{C}{D}{r}{1+\beta}+\beta r \right),
\end{array}
\end{eqnarray}

taking identity as the differential operator $\mathcal{B}$ we can build the differential equation

\begin{eqnarray}\label{eq:43}
\small
\left\{
\begin{array}{cl}
\mathcal{L}u = f,& \Omega,\\
\mathcal{B} u = g, &\partial \Omega,
\end{array}\right.=\left\{
\begin{array}{cl}
\mathcal{L}u = f,& \Omega,\\
 u = g, &\partial \Omega,
\end{array}\right.
\end{eqnarray}

where $ \Omega\subset \nset{R}^d$, in such a way that when $\beta \to 0$ and $d=2$  the equation \eqref{eq:43} takes the form of a Poisson's equation

\begin{eqnarray*}
\small
\left\{
\begin{array}{cl}
\nabla^2 u = f,& \Omega,\\
 u= g, &\partial \Omega,
\end{array}\right.
\end{eqnarray*}

then taking $d=2$ we get the following differential equation

\begin{eqnarray}\label{eq:48}
\small
\left\{
\begin{array}{cl}
\mathcal{L}u(x,y) = f(x,y),& \Omega,\\
 u(x,y) = g(x,y), &\partial \Omega,
\end{array}\right.
\end{eqnarray}

using the following distribution of interior nodes type Halton
together cartesian nodes near the boundary on the domain $\Omega_{0,1}$

\begin{figure}[!ht]
\includegraphics[width=0.45\textwidth]{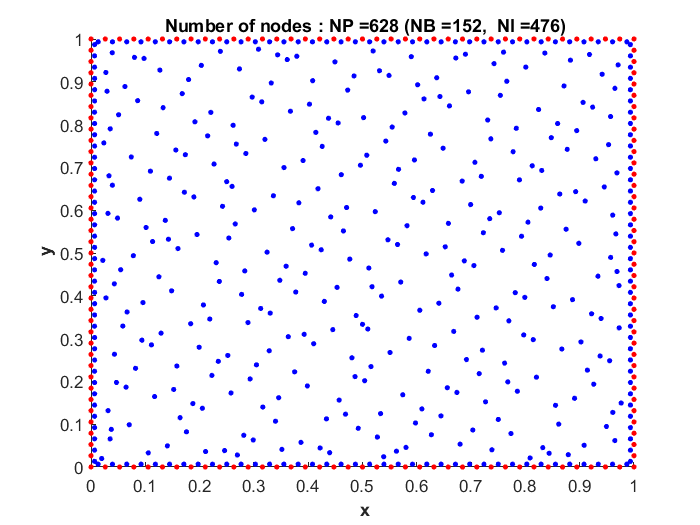}
\caption{Nodes used for the interpolation problem, where $NB$ and $NI$ are the boundary and interior nodes respectively
}\label{fig:02}
\end{figure}

Then to carry out the asymmetric collocation method is taken  $\alpha \in (-2,2)$, here is presented the option to use $ \alpha $ for fixed values although you can also use it looking for a value that minimizes the error.

Denoting by $\mathcal{L}\sigma_i=\mathcal{L}\sigma(x_i,y_i)$ and $f_i=f(x_i,y_i)$, then the error that we will use will be the root of the mean square error given by

\begin{eqnarray}\label{eq:49}
RMSE:=\sqrt{ \dfrac{1}{N_p}\sum_{i=1}^{N_p}\left(f_i-\mathcal{L}\sigma_i \right)^2}.
\end{eqnarray}

Taking $ \beta = -0.5 $ and the functions

\begin{eqnarray*}
\scriptsize
\begin{array}{c}
\begin{array}{ll}
f(x,y)=& \dfrac{2(108x - 36)^2(\cos(5.4y) + 1.25)}{(6(3x - 1)^2 + 6)^3}  - 
\dfrac{108\cos(5.4y) + 135}{(6(3x - 1)^2 + 6)^2},\\
&- \dfrac{29.16\cos(5.4y)}{6(3x - 1)^2 + 6},\\
g(x,y)=& \dfrac{\cos(5.4y)+1.25}{6(3x-1)^2+6} ,
\end{array}
\end{array}
\end{eqnarray*}

the asymmetrical collocation method is carried out

\begin{itemize}
\item[•] Using the radial function

\begin{eqnarray*}
\small{
\begin{array}{ll}
\Phi(\alpha, N,r)=&-2 b^{-N+\alpha}\ifr{0}{}{D}{r}{\alpha}r^{N+2}+4b^{1-N+\alpha}\ifr{0}{}{D}{r}{\alpha} r^{N+1}\\
&-2b^{2-N+\alpha}\ifr{0}{}{D}{r}{\alpha}r^N,
\end{array}
}
\end{eqnarray*}

to use the interpolant \eqref{eq:45} is taken

\begin{eqnarray*}
\small
\begin{array}{ll}
q&=\max\set{ord(\mathcal{L}),ord(\mathcal{B})}\\
&=\max\set{2+\beta,0}=1.5,
\end{array}
\end{eqnarray*} 

as $ \alpha \in (-2,2) $ is chosen

\begin{eqnarray*}
q+\alpha=3.5<N=3.55,
\end{eqnarray*}
 
and it is defined

\begin{eqnarray*}
\small
\begin{array}{l}
m=\max \set{\left\lceil \dfrac{N-(-2)+2}{2} \right\rceil,\left\lceil \dfrac{N-2+2}{2} \right\rceil }=4,
\end{array}
\end{eqnarray*}

obtaining the following results

\begin{eqnarray*}
\footnotesize
\begin{array}{ccc}
\alpha& RMSE &cond(G_M) \\ \hline
-1.9	&	9.4338346619063748e-2	&	7.6307799272592067	\\
 -1.8	&	8.2732749317909776e-2	&	5.3519676849454063	\\
 -1.7	&	8.8396716685719204e-2	&	8.5115507679137714	\\
 -1.6	&	8.1677337534535599e-2	&	4.0638967604070082	\\
 -1.5	&	8.7250610786736724e-2	&	6.0086096794810118	\\
 -1.4	&	1.8542238945053083e-1	&	3.9308669746402285	\\
 -1.3	&	3.7466868706221845e-1	&	6.1885465337305714	\\
 -1.0	&	4.2234223781422853e-1	&	7.5017560800341334	\\
 -0.7	&	6.0026543548922384e-1	&	5.5132119437845537	\\
 -0.6	&	8.8815582039528251e-1	&	8.1754571352817731	\\
  0.0	&	6.7896737435848153e-1	&	4.8990999804076649	\\
  0.1	&	2.0965746384065434e-1	&	7.1321575816746696	\\
  0.2	&	2.5755495145969060e-1	&	4.3454335303035379	\\
  0.3	&	3.1894520857777642e-1	&	8.0012752017703548	\\
  0.4	&	3.6725002710116900e-1	&	5.9760475954549683	\\
  0.5	&	4.2098721133896294e-1	&	6.1862716077593101	\\
  0.6	&	5.2079098187520356e-1	&	8.3771880377128802	\\
  0.7	&	6.7349862164498941e-1	&	3.7639598037087736					
\end{array}
\end{eqnarray*}

\newpage

\begin{figure}[!ht]
\includegraphics[width=0.5\textwidth]{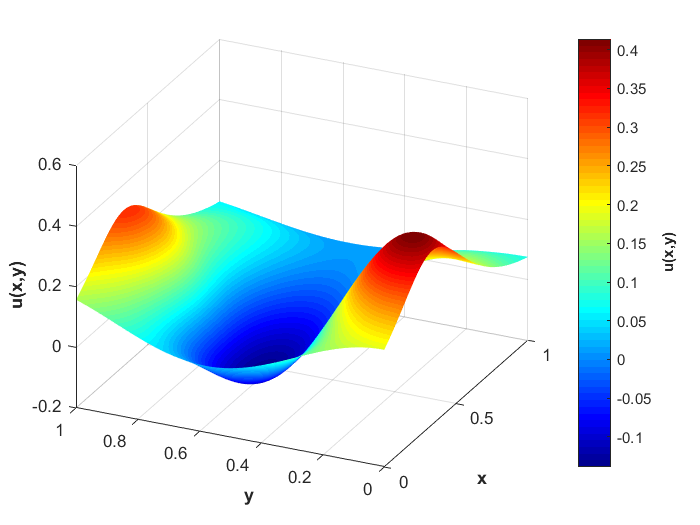}
\caption{Graph of the numerical solution (with minimal error) to the problem raised. }
\end{figure}

\end{itemize}

Taking $ \beta = 0.15 $ and the functions

\begin{eqnarray*}
\footnotesize
\begin{array}{c}
\begin{array}{ll}
f(x,y)=&-\dfrac{128}{35}\left(\sin(8(x+y))+\cos(8(x-y)) \right) ,\\
g(x,y)=& \dfrac{1}{35}\left(\sin(8(x+y))+\cos(8(x-y)) +4\right) ,
\end{array}
\end{array}
\end{eqnarray*}

the asymmetrical collocation method is carried out

\begin{itemize}
\item[•] Using the radial function

\begin{eqnarray*}
\small{
\begin{array}{ll}
\Phi(\alpha, N,r)=&-2 b^{-N+\alpha}\ifr{0}{}{D}{r}{\alpha}r^{N+2}+4b^{1-N+\alpha}\ifr{0}{}{D}{r}{\alpha} r^{N+1}\\
&-2b^{2-N+\alpha}\ifr{0}{}{D}{r}{\alpha}r^N,
\end{array}
}
\end{eqnarray*}

to use the interpolant \eqref{eq:45} is taken
 
\begin{eqnarray*}
\small
\begin{array}{ll}
q&=\max\set{ord(\mathcal{L}),ord(\mathcal{B})}\\
&=\max\set{2+\beta,0}=2.15,
\end{array}
\end{eqnarray*} 

as $ \alpha \in (-2,2) $ is chosen

\begin{eqnarray*}
q+\alpha=4.15<N=4.255,
\end{eqnarray*}
 
and it is defined

\begin{eqnarray*}
\small
\begin{array}{l}
m=\max \set{\left\lceil \dfrac{N-(-2)+2}{2} \right\rceil,\left\lceil \dfrac{N-2+2}{2} \right\rceil }=5,
\end{array}
\end{eqnarray*}

obtaining the following results

\begin{eqnarray*}
\footnotesize
\begin{array}{ccc}
\alpha& RMSE & cond(G_M)\\ \hline
 -1.6	&	6.5170655092245167e-1	&	4.7495389131006629	\\
 -1.5	&	1.7571501485532573e-1	&	4.0874067249340458	\\
 -1.4	&	1.2048058149601308e-1	&	7.1257745733005473	\\
 -1.3	&	1.2074660635701032e-1	&	9.1420630917154426	\\
 -1.2	&	1.2303027510252981e-1	&	7.7348875980860550	\\
 -1.1	&	1.4446689314716657e-1	&	7.6183854413772387	\\
 -1.0	&	1.8707079697143097e-1	&	4.7282991685218461	\\
 -0.9	&	2.8909651631485506e-1	&	6.5644036204148888	\\
 -0.8	&	6.5312833189177977e-1	&	5.4360430267446311	\\
 -0.7	&	8.3309014792233960e-1	&	7.9421114559403101	\\
 -0.6	&	3.6848119793984679e-1	&	8.3247898780765333	\\
 -0.1	&	4.3437390370092260e-1	&	4.5937327195633149	\\
  0.1	&	3.0172193192986441e-1	&	9.6741975397704021	\\
  0.2	&	3.4394425235516812e-1	&	4.5159570973161465	\\
  0.3	&	2.0982933894818576e-1	&	9.6547180568412188	\\
  0.4	&	2.4381512919793669e-1	&	5.1869902693095336	\\
  0.5	&	2.8971359130471058e-1	&	3.8030916671631547	\\
  0.6	&	3.4847114679571661e-1	&	8.2919349466597989	\\
  0.7	&	4.2320675674349861e-1	&	7.7574695544779200	\\
  0.8	&	5.1912850370692976e-1	&	5.9840232769986894	\\
  0.9	&	6.4660590409533081e-1	&	3.7594165306617882		
\end{array}
\end{eqnarray*}

\begin{figure}[!ht]
\includegraphics[width=0.5\textwidth]{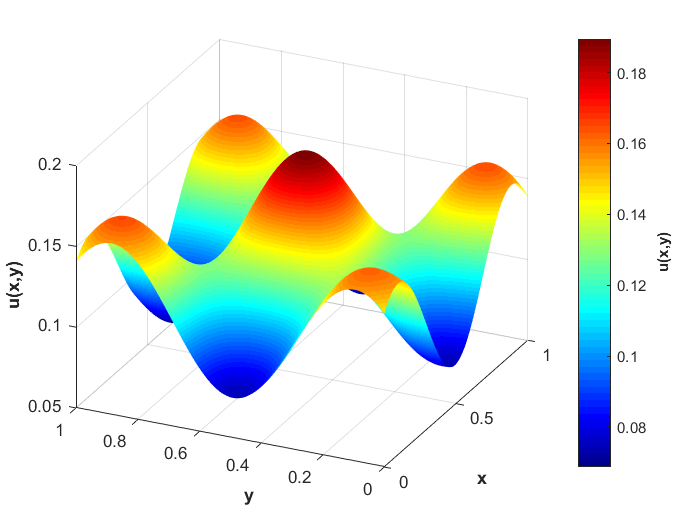}
\caption{Graph of the numerical solution (with minimal error) to the problem raised. }
\end{figure}

\end{itemize}

We can construct the next differential operator using the fractional derivative of Riemann-Liouville

\begin{eqnarray}
\small
\begin{array}{l}
\mathcal{L}=\left(\ifr{0}{}{D}{r}{2+\beta} + \dfrac{1}{r} \ifr{0}{}{D}{r}{1+\beta}+\beta r \right),
\end{array}
\end{eqnarray}

defined on a domain $\Omega \subset \nset{R}^d \setminus \set{0}$, taking identity as the differential operator $\mathcal{B}$ and  $d=2$ we get the following differential equation

\begin{eqnarray}
\small
\left\{
\begin{array}{cl}
\mathcal{L}u(x,y) = f(x,y),& \Omega,\\
 u(x,y) = g(x,y), &\partial \Omega,
\end{array}\right.
\end{eqnarray}

using the following distribution of interior nodes type Halton
together cartesian nodes near the boundary on the domain $\Omega_{0.28,1.48}$

\begin{figure}[!ht]
\includegraphics[width=0.45\textwidth]{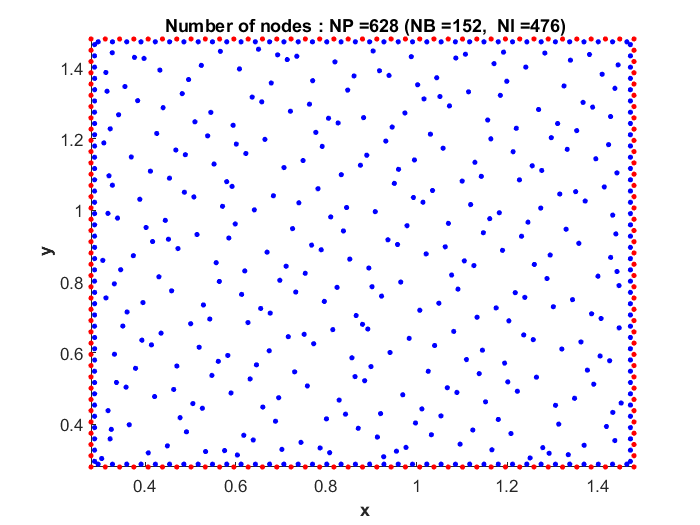}
\caption{Nodes used for the interpolation problem, where $NB$ and $NI$ are the boundary and interior nodes respectively
}\label{fig:03}
\end{figure}

then taking $ \beta = -2.5 $ and the functions

\begin{eqnarray*}
\footnotesize
\begin{array}{c}
\begin{array}{ll}
f(x,y)=&-\dfrac{128}{35}\left(\sin(8(x+y))+\cos(8(x-y)) \right) ,\\
g(x,y)=& \dfrac{1}{35}\left(\sin(8(x+y))+\cos(8(x-y)) +4\right) ,
\end{array}
\end{array}
\end{eqnarray*}

the asymmetrical collocation method is carried out

\begin{itemize}
\item[•] Using the radial function

\begin{eqnarray*}
\small{
\begin{array}{ll}
\Phi(\alpha, N,r)=&-2 b^{-N+\alpha}\ifr{0}{}{D}{r}{\alpha}r^{N+2}+4b^{1-N+\alpha}\ifr{0}{}{D}{r}{\alpha} r^{N+1}\\
&-2b^{2-N+\alpha}\ifr{0}{}{D}{r}{\alpha}r^N,
\end{array}
}
\end{eqnarray*}

to use the interpolant \eqref{eq:45} is taken
 
\begin{eqnarray*}
\small
\begin{array}{ll}
q&=\max\set{ord(\mathcal{L}),ord(\mathcal{B})}\\
&=\max\set{2+\beta,0}=0,
\end{array}
\end{eqnarray*} 

as $ \alpha \in (-2,2) $ is chosen

\begin{eqnarray*}
\alpha=2<N=2.25,
\end{eqnarray*}
 
and it is defined

\begin{eqnarray*}
\small
\begin{array}{l}
m=\max \set{\left\lceil \dfrac{N-(-2)+2}{2} \right\rceil,\left\lceil \dfrac{N-2+2}{2} \right\rceil }=4,
\end{array}
\end{eqnarray*}
obtaining the following results

\begin{eqnarray*}
\footnotesize
\begin{array}{ccc}
\alpha& RMSE & cond(G_M)\\ \hline
-1.3	&	8.1624992247252071e-1	&	4.8877253453360776	\\
 -1.2	&	8.7150356043661481e-1	&	7.8879861055952984	\\
 -1.0	&	2.3947625681238222e-1	&	4.6563142060453906	\\
 -0.9	&	6.3813942503541832e-1	&	6.7415090485234304	\\
 -0.8	&	2.4245563889737556e-1	&	9.0095667622185402	\\
 -0.7	&	7.8273416767854806e-2	&	4.6153771286769230	\\
 -0.6	&	6.8520649294367561e-2	&	4.8689277903455732	\\
 -0.5	&	1.2800563138668500e-1	&	4.0912301014848360	\\
 -0.4	&	4.3817929756525886e-1	&	7.9578762827989884	\\
 -0.3	&	4.7076551391474242e-1	&	6.3858782448176070	\\
 -0.1	&	5.1458635004984754e-1	&	6.0577940905942782	\\
  0.0	&	7.7818081533588657e-1	&	3.8807153590997578	\\
  0.1	&	6.7401786461338986e-1	&	5.2586169156509017	\\
  0.3	&	4.5597146882611150e-1	&	4.3704175065445598	\\
  0.5	&	2.8305679271014683e-1	&	6.1044635630032165	\\
  0.6	&	6.0966123979559073e-1	&	8.8043587420730951	\\
  0.7	&	1.8786815492663200e-1	&	3.9477776389878581	\\
  0.8	&	4.0552430446664012e-1	&	6.7368767950043269	\\
  0.9	&	2.1706046218673369e-1	&	7.0190727229648013	\\
  1.0	&	4.7725454402349632e-1	&	8.4827242519603523	\\
  1.1	&	8.0376048230088470e-1	&	4.2109594173989944			
\end{array}
\end{eqnarray*}

\begin{figure}[!ht]
\includegraphics[width=0.5\textwidth]{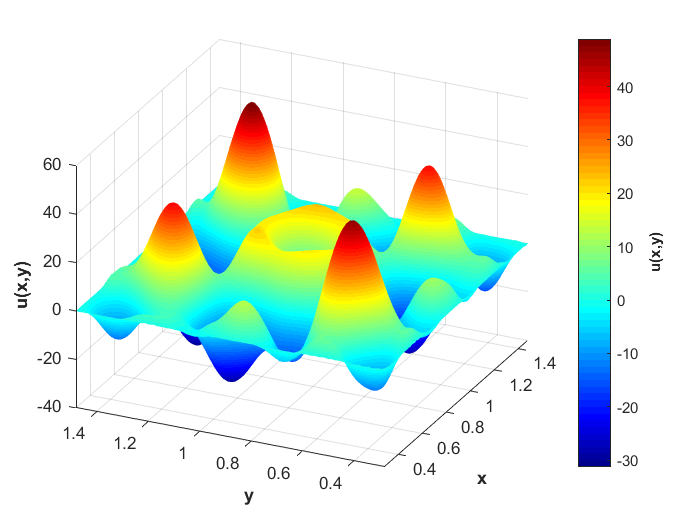}
\caption{Graph of the numerical solution (with minimal error) to the problem raised. }
\end{figure}

\end{itemize}

Although examples are presented in $ \nset{R}^2 $ the radial functions presented in this document are valid for $\nset{R}^d$, to implement them we can assume that have a domain $\Omega\subset \Omega_d \subset\nset{R}^d$, where

\begin{eqnarray*}
\small
\begin{array}{l}
\Omega_d= [a_1,b_1]\times [a_2,b_2]\times \cdots \times [a_d,b_d],
\end{array}
\end{eqnarray*}

then it is enough to define

\begin{eqnarray*}
\small
\begin{array}{l}
b=\max\set{ \set{b_i}_{i=1}^d}.
\end{array}
\end{eqnarray*}

Although the definition of fractional derivative of Riemann-Liouville was used for the functions constructed in the previous sections, in general, any other definition of fractional derivative can be used as long as this definition is used at par with the fractional integral of Riemann-Liouville.

\bibliographystyle{unsrt}
\nocite{torres2019fractional}
\nocite{bramb17}
\nocite{fernando2017fractional}
\nocite{brambila2018fractional}
\nocite{martinez2017application2}
\nocite{martinez2017applications1}
\bibliography{Biblio}

\end{document}